# OPTIMAL RATES OF CONVERGENCE FOR ESTIMATING THE NULL DENSITY AND PROPORTION OF NONNULL EFFECTS IN LARGE-SCALE MULTIPLE TESTING

BY T. TONY CAI[1] AND JIASHUN JIN[2]

*University of Pennsylvania and Carnegie Mellon University*

An important estimation problem that is closely related to large-scale multiple testing is that of estimating the null density and the proportion of nonnull effects. A few estimators have been introduced in the literature; however, several important problems, including the evaluation of the minimax rate of convergence and the construction of rate-optimal estimators, remain open.

In this paper, we consider optimal estimation of the null density and the proportion of nonnull effects. Both minimax lower and upper bounds are derived. The lower bound is established by a two-point testing argument, where at the core is the novel construction of two least favorable marginal densities $f_1$ and $f_2$. The density $f_1$ is heavy tailed both in the spatial and frequency domains and $f_2$ is a perturbation of $f_1$ such that the characteristic functions associated with $f_1$ and $f_2$ match each other in low frequencies. The minimax upper bound is obtained by constructing estimators which rely on the empirical characteristic function and Fourier analysis. The estimator is shown to be minimax rate optimal.

Compared to existing methods in the literature, the proposed procedure not only provides more precise estimates of the null density and the proportion of the nonnull effects, but also yields more accurate results when used inside some multiple testing procedures which aim at controlling the False Discovery Rate (FDR). The procedure is easy to implement and numerical results are given.

**1. Introduction.** Large-scale multiple testing is an important area in modern statistics with a wide range of applications including DNA microar-

Received November 2008; revised February 2009.
[1]Supported in part by NSF Grant DMS-06-04954.
[2]Supported in part by NSF Grants DMS-05-05423 and DMS-09-08613.

ray studies, functional Magnetic Resonance Imaging analyses (fMRI) and astronomical surveys. Since the seminal paper by Benjamini and Hochberg (1995) on false discovery rate (FDR) control, research in this area has been very active. See, for example, Efron et al. (2001), Storey (2002), Genovese and Wasserman (2004), van der Laan, Dudoit and Pollard (2004) and Sun and Cai (2007). Properties of FDR-controlling procedures have been studied, for example, in Finner, Dickhaus and Roters (2009) and Neuvial (2008). See also Abramovich et al. (2006) and Donoho and Jin (2006) for estimation using a multiple testing approach.

In large-scale multiple testing, one tests *simultaneously* a large number of null hypotheses

$$H_1, H_2, \ldots, H_n. \tag{1.1}$$

Frequently, associated with each hypothesis $H_j$ is a test statistic $X_j$, which can be a $z$-score, a $p$-value, a summary statistic, etc., depending on the situation. The goal is to use the test statistics to determine which hypotheses are true and which are false. We call $X_j$ a *null effect* if $H_j$ is true and a *nonnull effect* otherwise.

A commonly used and effective framework for large-scale multiple testing is the so-called two-group random mixture model which assumes that each hypothesis has a given probability of being true and the test statistics are generated from a mixture of two densities; see, for example, Efron et al. (2001), Newton et al. (2001), Storey (2002) and Sun and Cai (2007). In detail, let $\theta = (\theta_1, \ldots, \theta_n)$ be independent Bernoulli($\varepsilon$) variables, where $\varepsilon \in (0,1)$ and $\theta_j = 0$ indicates that the null hypothesis $H_j$ is true and $\theta_j = 1$ otherwise. When $\theta_j = 0$, $X_j$ is generated from a density $f^{\text{null}}(x)$. When $\theta_j = 1$, $X_j$ is generated from another (alternative) density $f^{\text{alt}}(x)$. Marginally, $X_j$ obeys the following two-group random mixture model:

$$X_j \overset{\text{i.i.d.}}{\sim} (1-\varepsilon)f^{\text{null}} + \varepsilon f^{\text{alt}} \equiv f, \qquad j=1,\ldots,n, \tag{1.2}$$

where $f^{\text{null}}$, $f^{\text{alt}}$ and $\varepsilon$ are called the *null density*, *nonnull density* and *proportion of nonnull effects*, respectively.

An important estimation problem that is closely related to multiple testing is that of estimating $f^{\text{null}}$, $\varepsilon$ and $f$. In fact, many commonly used multiple testing procedures require good estimators of some or all of these three quantities. See Benjamini and Hochberg (2000), Efron et al. (2001), Storey (2002), Genovese and Wasserman (2004), Benjamini, Krieger and Yekutieli (2006), Blanchard and Roquain (2007) and Sun and Cai (2007). For example, in an empirical Bayes framework, Efron et al. (2001) introduced the local false discovery rate (Lfdr) which is defined as

$$\text{Lfdr}(x) = \frac{(1-\varepsilon)f^{\text{null}}(x)}{f(x)}. \tag{1.3}$$

Lfdr has a useful Bayesian interpretation as the a posteriori probability of a hypothesis being in the null group given the value of the test statistic. See also Müller et al. (2004). Sun and Cai (2007) considered the multiple testing problem from a compound decision theoretical point of view and showed that the Lfdr is a fundamental quantity which can be used directly for the optimal FDR control. Calculating the Lfdr clearly requires the knowledge of $\varepsilon$, $f^{\text{null}}$ and $f$. In real applications, the proportion $\varepsilon$ and the marginal density $f$ are unknown and thus need to be estimated from the data. The null density $f^{\text{null}}$ is more subtle. In many studies the null distribution is assumed to be known and can be used directly for multiple testing. However, somewhat surprisingly, Efron (2004) demonstrated convincingly that in some applications such as the analysis of microarray data on breast cancer and human immunodeficiency virus (HIV) the true null distribution of the test statistic can be quite different from the theoretical null, and possible causes for such a phenomenon include but are not limited to unobserved covariates, correlations across different arrays and different genes. It is further illustrated in Jin and Cai (2007) that two seemingly close choices of the null distribution can lead to substantially different testing results. Hence, a careful study on how to estimate the null distribution is also indispensable.

In the present paper we study the problem of optimal estimation of the null density $f^{\text{null}}$ and the proportion $\varepsilon$. We should mention that estimating the marginal density $f$ is a standard density estimation problem and is well understood. See, for example, Silverman (1986). Several methods for estimating the null density $f^{\text{null}}$ and the proportion $\varepsilon$ have been introduced in the literature. See Efron (2004, 2008) and Jin and Cai (2007) for estimating $f^{\text{null}}$ and $\varepsilon$, and see Genovese and Wasserman (2004), Meinshausen and Rice (2006), Cai, Jin and Low (2007), Jin (2008) and Celisse and Robin (2008) for estimating $\varepsilon$ [also see Storey (2002), Efron et al. (2001), Swanepoel (1999)]. Unfortunately, despite the encouraging progress in these works, the optimality of the estimators is largely unknown [it is, however, not hard to show that some of these estimators are generally inconsistent in the nonsparse case; see, e.g., Jin and Cai (2007)]. It is hence of significant interest to understand how well $f^{\text{null}}$ and $\varepsilon$ can be estimated and to what extend improving the estimation accuracy of $f^{\text{null}}$ and $\varepsilon$ can help to enhance the performance of leading contemporary multiple testing procedures [including but not limited to those by Benjamini and Hochberg (1995), Efron et al. (2001) and Sun and Cai (2007)]. Multiple testing procedures that adapt to $\varepsilon$, without estimating it directly, have also been proposed recently in Blanchard and Roquain (2007) and Finner, Dickhaus and Roters (2009).

In this paper, we focus on the Gaussian mixture model as in Efron (2004). We model $f^{\text{null}}$ as Gaussian, but both the mean and the variance are un-known and need to be estimated:

$$(1.4) \qquad f^{\text{null}}(x) = \frac{1}{\sigma_0}\phi\left(\frac{x-u_0}{\sigma_0}\right), \qquad \phi\text{: density of } N(0,1).$$

We shall use the terminology in Efron (2004) by calling $\sigma_0^2$ the *null variance parameter*, $u_0$ the *null mean parameter*, and together the *null parameters*. The Gaussian model for $f^{\text{null}}$ is somewhat idealized, but it is a reasonable choice. On one hand, assuming $f^{\text{null}}$ as Gaussian helps to re-normalize the null distribution and is therefore a good starting point in large-scale multiple testing. On the other hand, allowing $f^{\text{null}}$ to be in a much broader class will lead to identifiability problems. The nonnull distribution $f^{\text{alt}}$ is modeled by a Gaussian location-scale mixture,

$$(1.5) \qquad f^{\text{alt}}(x) = \int\int \frac{1}{\sigma}\phi\left(\frac{x-u}{\sigma}\right) dH(u,\sigma),$$

where $H$ is called the *mixing distribution*. Additional to the mathematical tractability that it offers, model (1.5) also offers great flexibility. For example, it is well known that under the $L^1$-metric, the set of Gaussian mixing densities of the form in (1.5) is dense in the set of all density functions. Also, model (1.5) is able to capture the essence of many application examples. See Jin (2008) for an example on the analysis of gene microarray data on breast cancer and an example on the study of the abundance of the Kuiper Belt Objects.

We consider the asymptotic minimax estimation problem and address several inter-connected questions: what are the optimal rates of convergence? what are the best estimation tools, and where do the difficulties of the estimation problem come from? Our analysis reveals that the optimal rates of convergence for estimating the proportion and the null parameters depend on the smoothness of $H(u,\sigma)$ (more specifically, the conditional density of $u$ given $\sigma$ associated with $H$). For an intuitive explanation, we note that $f^{\text{null}}$ and $f^{\text{alt}}$ are the convolution of the standard Gaussian with the point mass concentrated at $(u_0,\sigma_0)$ and $H$, respectively. Therefore, the smoother $H$ is, the more "different" it is from a point mass, and the less similar that $f^{\text{null}}$ and $f^{\text{alt}}$ are. Consequently, it is easier to separate one from the other, and hence a faster convergence rate in estimating the proportion and the null parameters.

Since the smoothness of a density can be conveniently characterized by the tail behavior of its characteristic function, this suggests that frequency domain techniques can be naturally used for studying the optimal rate of convergence. Along this line, we first derive a minimax lower bound by a careful analysis of the tail behavior of the characteristic functions and by a two-point testing technique. We then establish the upper bound by constructing estimators with the risks converging to zero at the same rate as that of the lower bound—such estimators are then rate optimal. The procedures are closely related to our recent work Jin and Cai (2007) and Jin (2008) which to the best of our knowledge are the only frequency-domain-based approach to estimating the null parameters and the proportion of nonnull effects. We should emphasize that the upper bound does not follow trivially from that in Jin and Cai (2007) and Jin (2008). For example, it is seen that the procedure for estimating the proportion proposed in Jin and Cai (2007) and Jin (2008) is not optimal, and careful modifications are needed to make it optimal. Also, to prove the optimality of the procedures here, we need much more delicate analysis than that in Jin and Cai (2007) and Jin (2008), where the scope of the study is limited to the consistency of the procedures.

In addition to the asymptotic analysis, we also investigate the finite $n$ performance of the estimators using simulated data. The proposed procedures are easy to implement. The goal for the simulation study is two-fold: how accurate the parameters are estimated and how the errors in the point estimation affect the results of the subsequent multiple testing. The numerical study shows that our estimators enjoy superior performance both in parameter estimation (measured by mean squared errors) and in the subsequent multiple testing. Our estimator of the proportion performs well uniformly in all the cases in comparison to the estimators proposed in Storey (2002) and Efron (2004). In particular, it is robust under many different choices of nonnull distribution and sparsity level. The multiple testing results are generally sensitive to the changes in the null parameters as well as the proportion. In our numerical study, we compare the performance of our estimators with those of Storey (2002) and Efron (2004) using two specific multiple testing procedures, the adaptive $p$-value based procedure of Benjamini and Hochberg (2000) which requires estimation of the proportion $\varepsilon$, and the AdaptZ procedure of Sun and Cai (2007) which requires estimation of $\varepsilon$, $f$ and $f^{\text{null}}$. The simulation study shows that our estimators yield the most accurate multiple testing results in both cases in comparison to the other two estimators.

The paper is organized as follows. In Section 2, after basic notation and definitions are introduced, we consider the minimax lower bound for estimating the null parameters. We then derive the minimax rates of convergence by showing that the lower bound is in fact sharp. This is accomplished by constructing rate-optimal estimators using the empirical characteristic functions. Section 3 studies the minimax estimation of the proportion. We first consider the simpler case where the null parameters are given and then extend the result to the case where the null parameters are unknown. Section 4 investigates the numerical performance of our procedure by a simulation study. Section 5 discusses possible extensions of our work and its connections with the nonparametric deconvolution problem. The proofs of the main results are given in Section 6 and the Appendix contains the proofs of the technical lemmas that are used to prove the main results.

**2. Estimating the null parameters: Minimax risk and rate optimal estimators.** In this section, we study the minimax risks for estimating the null parameters. The minimax lower bounds are established by a two-point testing argument in Section 2.1. At the core of the argument is the construction of two underlying densities whose corresponding null parameters are different but whose characteristic functions match with each other in low frequencies. We then derive the minimax upper bounds by constructing and studying rate optimal estimators in Section 2.2.

Return to the Gaussian mixture model

$$(2.1) \quad X_j \stackrel{\text{i.i.d.}}{\sim} (1-\varepsilon)\frac{1}{\sigma_0}\phi\left(\frac{x-u_0}{\sigma_0}\right) + \varepsilon \int \frac{1}{\sigma}\phi\left(\frac{x-u}{\sigma}\right) dH(u,\sigma) \equiv f(x).$$

For any mixing distribution $H(u,\sigma)$ under consideration, let $H(\sigma)$ be the marginal distribution of $\sigma$ and let $H(u|\sigma)$ be the conditional distribution of $u$ given $\sigma$.

DEFINITION 2.1. We call a density $f$ eligible if it has the form as in (2.1) where $H(u,\sigma)$ satisfies that $H(\sigma)$ is supported on $[\sigma_0,\infty)$ and that $H(u|\sigma)$ has a density $h(u|\sigma)$ for any $\sigma \geq \sigma_0$. We denote the set of all eligible $f$ by $\mathcal{F}$.

Two examples for eligible $f$ are (1). $H(\sigma)$ is supported in $[\sigma_0 + \delta, \infty)$ for some constant $\delta > 0$, and (2). $H(\sigma)$ is the point mass at $\sigma_0$, and $H(u|\sigma_0)$ has a density.

In this paper, we focus on eligible $f$, so that the null parameters and the proportion of nonnull effects are both identifiable. See Jin and Cai (2007) for more discussion on identifiability.

We shall define the parameter space of $f$ for the minimax theory. First, we suppose that for some fixed constant $q > 0$ and $A > a > 0$,

$$(2.2) \quad \sigma_0 \geq a, \qquad \int |x|^q f(x)\, dx \leq A^q,$$

so that $\sigma_0^2$ and $u_0$ are uniformly bounded across the whole parameter space. Second, fix $\alpha > 0$. We assume

$$(2.3) \quad \overline{\lim_{t \to \infty}} \sup_{\sigma \geq \sigma_0} \{|t|^\alpha |\hat{h}(t|\sigma)|\} \leq A, \qquad \overline{\lim_{t \to \infty}} \sup_{\sigma \geq \sigma_0} \{|t|^{\alpha+1}|\tilde{h}'(t|\sigma)|\} \leq A,$$

where $h(u|\sigma)$ is the aforementioned conditional density, $\hat{h}(t|\sigma)$ is the corresponding characteristic function and

$$(2.4) \quad \tilde{h}(t|\sigma) = \tilde{h}(t|\sigma; u_0) = \int e^{itu} h(u+u_0|\sigma)\, du.$$

Roughly speaking, (2.3) requires $h(u|\sigma)$ to be sufficiently smooth so that $\hat{h}(t|\sigma)$ decays at a rate not slower than that of $|t|^{-\alpha}$. We shall see below that the minimax risk depends on the smoothness parameter $\alpha$. Note that in (2.2) and (2.3), different constants $A$ can be used in different places. However, this does not change the minimax rate of convergence, so we use the same $A$ for simplicity.

Last, we calibrate the proportion $\varepsilon$. In the literature, the proportion is a well-known measure for *sparsity*; see, for example, Abramovich et al. (2006) and Jin and Cai (2007). In this paper, we focus on the *moderately sparse case* where the proportion $\varepsilon = \varepsilon_n$ can be small but not smaller than $1/\sqrt{n}$. The case $\varepsilon_n \ll 1/\sqrt{n}$ is called the *very sparse case* and has been proven to be much more challenging for statistical inference; see Donoho and Jin (2004) and Cai, Jin and Low (2007) for detailed discussion. In light of this, we suppose that for some fixed parameters $\varepsilon_0 \in (0,1)$ and $\beta \in [0, 1/2)$,

$$(2.5) \qquad \varepsilon_n \leq \eta_n \qquad \text{where } \eta_n = \eta_n(\varepsilon_0, \beta) \equiv \varepsilon_0 n^{-\beta}.$$

Note that $\eta_n = \varepsilon_0$ when $\beta = 0$. For this reason, we require $\varepsilon_0 < 1$ so that the null component will not be vanishingly small.

In summary, the parameter space we consider for the minimax risk is

$$(2.6) \quad \begin{aligned} \mathcal{F}_0 &= \mathcal{F}_0(\alpha, \beta, \varepsilon_0, q, a, A; n) \\ &= \{f \in \mathcal{F} \text{ and satisfies } (2.2), (2.3) \text{ and } (2.5)\}. \end{aligned}$$

We measure the performance of an estimator for the null parameters by mean squared errors, and measure the level of difficulty for the problem of estimating the null parameters $\sigma_0^2$ and $u_0$ by the minimax risks defined, respectively, by

$$R_n^\sigma = R_n^\sigma(\mathcal{F}_0(\alpha, \beta, \varepsilon_0, q, a, A; n)) = \inf_{\hat{\sigma}_0^2} \left\{ \sup_{\mathcal{F}_0(\alpha, \beta, \varepsilon_0, q, a, A; n)} E[\hat{\sigma}^2 - \sigma_0^2]^2 \right\}$$

and

$$R_n^u = R_n^u(\mathcal{F}_0(\alpha, \beta, \varepsilon_0, q, a, A; n)) = \inf_{\hat{u}_0} \left\{ \sup_{\mathcal{F}_0(\alpha, \beta, \varepsilon_0, q, a, A; n)} E[\hat{u}_0 - u_0]^2 \right\}.$$

2.1. *Lower bound for the minimax risk.* In this section, we establish the lower bound for the minimax risk of estimating $\sigma_0^2$ and $u_0$. As the discussions are similar, we shall focus on that for $\sigma_0^2$. We use the well-known two-point testing argument to show the lower bound [see, e.g., Ibragimov, Nemirovskii and Khas'minskii (1986) and Donoho and Liu (1991)], where the key is to construct two density functions in $\mathcal{F}_0$—$f_1(x)$ and $f_2(x)$—such that the null variance parameters associated with them differ by a small amount, say $\delta_n$, but two densities are indistinguishable in the sense that their $\chi^2$-distance

$$(2.7) \qquad d(f_1, f_2) \equiv \int \frac{(f_2(x) - f_1(x))^2}{f_1(x)} \, dx$$

is of a smaller order than that of $1/n$. In fact, once such densities $f_1$ and $f_2$ are constructed, then there is a constant $C > 0$ such that

$$R_n^\sigma \geq C\delta_n^2 \tag{2.8}$$

and $C\delta_n^2$ is a lower bound for the minimax risk; see Ibragimov, Nemirovskii and Khas'minskii (1986) and Donoho and Liu (1991) for details.

To this end, let

$$a_n^2 = a^2 + \delta_n,$$

where $\delta_n > 0$ to be determined. Our construction of $f_1$ and $f_2$ has the form of

$$f_1(x) = (1 - \eta_n)\frac{1}{a}\phi\left(\frac{x}{a}\right) + \eta_n \int \frac{1}{a}\phi\left(\frac{x-u}{a}\right)h_1(u)\,du, \tag{2.9}$$

$$f_2(x) = (1 - \eta_n)\frac{1}{a_n}\phi\left(\frac{x}{a_n}\right) + \eta_n \int \frac{1}{a_n}\phi\left(\frac{x-u}{a_n}\right)h_2(u)\,du, \tag{2.10}$$

where $a$ and $\eta_n$ are as in the definition of $\mathcal{F}_0(\alpha, \beta, \varepsilon_0, q, a, A; n)$, $h_1(u)$ and $h_2(u)$ are two density functions to be determined (note that the null variance parameters associated with $f_1$ and $f_2$ differ by an amount of $\delta_n$). There are two key elements in our construction. First, the characteristic functions of $f_1$ and $f_2$ match with each other in low frequencies, that is, for a constant $\tau = \tau_n$ to be determined,

$$\hat{f}_1(t) = \hat{f}_2(t) \qquad \forall |t| \leq \tau_n. \tag{2.11}$$

Second, $f_1$ is heavy-tailed in the spatial domain,

$$f_1(x) \geq C\eta_n(1 + |x|)^{-k} \qquad \forall x, \tag{2.12}$$

where $k > 0$ is an integer to be determined. Below, we first show that the $\chi^2$-distance between $f_1$ and $f_2$ equals to $o(1/n)$ if we take the $\tau_n$ in (2.11) to be

$$\tau_n = \frac{1}{a}\sqrt{3\log n}. \tag{2.13}$$

We then sketch how to construct $f_1$ and $f_2$ to satisfy (2.11) and (2.12), and discuss how large $\delta_n$ could be so that such a construction is possible. We conclude this subsection with the statement for the minimax lower bounds. To focus on the main ideas, we try to be simple and heuristic in this section and leave proof details to Section 6.

We now begin by investigating the $\chi^2$-distance. First, the heavy-tailed property of $f_1$ largely simplifies the calculation of the $\chi^2$-distance. In fact, by (2.12) and the well-known Parseval formula [Mallat (1998)], the $\chi^2$-distance is proportional to the $L^2$-distance in the spatial domain, and so the $L^2$-distance in the frequency domain,

$$d(f_1, f_2) \le C \log^{k/2}(n) \eta_n^{-1} \int (f_2(x) - f_1(x))^2 \, dx$$
$$= C \log^{k/2}(n) \eta_n^{-1} \int (\hat{f}_1(t) - \hat{f}_2(t))^2 \, dt.$$

See Section 6 for the proof. Moreover, since that $\hat{f}_1$ and $\hat{f}_2$ match each other in low frequencies, and that $|\hat{f}_j(t)| \le Ce^{-a^2 t^2/2}$ for $j = 1, 2$,

$$\int (\hat{f}_1(t) - \hat{f}_2(t))^2 \, dt = \int_{|t| \ge \tau_n} (\hat{f}_1(t) - \hat{f}_2(t))^2 \, dt \le C \int_{|t| \ge \tau_n} e^{-a^2 t^2/2} \, dt.$$

Putting these together,

$$(2.14) \quad d(f_1, f_2) \le C \log^{k/2}(n) \eta_n^{-1} e^{-a^2 \tau_n^2/2} = C \eta_n^{-1} \log^{k/2}(n) n^{-3/2}.$$

Since $\eta_n \gg 1/\sqrt{n}$, this show that the $\chi^2$-distance $d(f_1, f_2) = o(1/n)$.

Next, we sketch the idea for constructing $f_1$ and $f_2$. Consider $f_1$ first. We construct $h_1$ as a perturbation of the standard normal density,

$$(2.15) \quad h_1(u) = \phi(u) + \vartheta_0 w_1(u).$$

The key is to show that for an appropriate constant $\vartheta_0 > 0$ and a function $w_1$, $h_1$ is indeed a density function, and $f_1$ satisfies the heavy-tailed requirement (2.12). Let $k$ be an even number, we construct $w_1(u)$ through its characteristic function as follows: $\hat{w}_1(t) = \frac{(-1)^{k/2}\pi}{(k-1)!}|t|^{k-1}$ in the vicinity of 0, $\hat{w}_1(t) = |t|^{-\alpha}$ for large $|t|$, and is smooth in between [details are given later in (6.1)]. By elementary Fourier analysis, first, we note that $\int w_1(u) \, du = \hat{w}_1(0) = 0$. Second, we note that the tail behavior of $w_1$ is determined by the only singular point of $\hat{w}_1$ (which is $t = 0$); in fact, by repeatedly using integration by parts, we have that for large $u$, $w_1(u) \sim |u|^{-k}$, that is,

$$(2.16) \quad \lim_{|u| \to \infty} w_1(u)|u|^k = 1.$$

We shall see that, first, (2.16) implies the heavy-tailed property of $f_1$, and second, (2.16) ensures that $w_1(u)$ is positive for sufficiently large $u$, so $h_1$ is a density function for an appropriately small $\vartheta_0 > 0$. Additionally, we will justify later that $f_1$ belongs to $\mathcal{F}_0$. Therefore, $f_1$ constructed this way meets all the desired requirements.

Now consider $f_2$. Similarly, we construct $h_2$ as a perturbation of a normal density,

$$(2.17) \quad h_2(u) = \frac{1}{\sqrt{1-\delta_n}} \phi\left(\frac{u}{\sqrt{1-\delta_n}}\right) + \vartheta_0 w_2(u)$$

and the key is to construct $w_2$ so that $\hat{f}_1$ and $\hat{f}_2$ match in low frequencies. Note that

$$\hat{f}_1(t) = \eta_n e^{-(a^2+1)t^2/2} + e^{-a^2 t^2/2}[(1-\eta_n) + \vartheta_0 \eta_n \hat{w}_1(t)]$$

and

$$\hat{f}_2(t) = \eta_n e^{-(a^2+1)t^2/2} + e^{-a_n^2 t^2/2}[(1-\eta_n) + \vartheta_0 \eta_n \hat{w}_2(t)].$$

By direct calculations, in order for $\hat{f}_1$ and $\hat{f}_2$ to match in low frequencies, it is necessary that

(2.18)
$$\hat{w}_2(t) = \tilde{w}(t)$$

for all $|t| \leq \tau_n$ where $\tilde{w}(t) \equiv e^{\delta_n t^2/2} \hat{w}_1(t) + \frac{1}{\vartheta_0} \frac{1-\eta_n}{\eta_n}[e^{\delta_n t^2/2} - 1].$

In light of this, we construct $w_2$ through its characteristic function as follows: $\hat{w}_2(t) = \tilde{w}(t)$ for $|t| \leq \tau_n$, $\hat{w}_2(t) = 0$ for $|t| > \tau_n + 1$, and is smooth in between. Figure 1 illustrates the construction of $\hat{w}_1$ and $\hat{w}_2$; see details therein.

We now investigate what is the largest $\delta_n$ so that $f_2$ constructed this way belongs to $\mathcal{F}_0$. By the definition of $\mathcal{F}_0$, it is necessary that $|\hat{h}_2(t)| \leq A|t|^{-\alpha}$ for all $t$, and especially that $|\hat{h}_2(\tau_n)| \leq A\tau_n^{-\alpha}$. Recall that $\hat{w}_1(\tau_n) = \vartheta_0 \tau_n^{-\alpha}$, we have

$$\hat{h}_2(\tau_n) = e^{\delta_n t^2/2} \hat{w}_1(\tau_n) + \frac{1}{\vartheta_0} \frac{1-\eta_n}{\eta_n}[e^{\delta_n \tau_n^2/2} - 1] \sim O\left(\tau_n^{-\alpha} + \frac{\delta_n}{\vartheta_0 \eta_n} \tau_n^2\right).$$

Together, these require that

$$\delta_n \leq C \eta_n \tau_n^{-(\alpha+2)}.$$

In light of this, we calibrate $\delta_n$ as

(2.19)
$$\delta_n = \theta_0 \vartheta_0 \eta_n \tau_n^{-(\alpha+2)},$$

where $\theta_0 > 0$ is a constant to be determined. Interestingly, it turns out that for an appropriately small $\theta_0$, $w_2$ constructed in this way ensures that $h_2$ is a density function and that $f_2$ lives $\mathcal{F}_0$ (see Section 6). Therefore, the largest possible $\delta_n$ is of the order of $O(\eta_n \tau_n^{-(\alpha+2)})$.

We are now ready to state the minimax lower bounds. Let $M_q$ be the $q$th moment of the standard normal [i.e., $M_q = E|X|^q$ with $X \sim N(0,1)$], the following theorem is proved in Section 6.

THEOREM 2.1. *Fix $\alpha > 2$, $\beta \in [0, 1/2)$, $\varepsilon_0 \in (0,1)$, $q > 0$, $a > 0$ and $A > \sqrt{a^2+1} M_q^{1/q}$. There is a constant $C > 0$ which depends on $\alpha, \beta, \varepsilon_0, q, a$ and $A$ such that,*

$$\lim_{n \to \infty} n^{2\beta} \cdot (\log n)^{(\alpha+2)} \cdot R_n^\sigma(\mathcal{F}_0(\alpha, \beta, \varepsilon_0, q, a, A; n)) \geq C$$

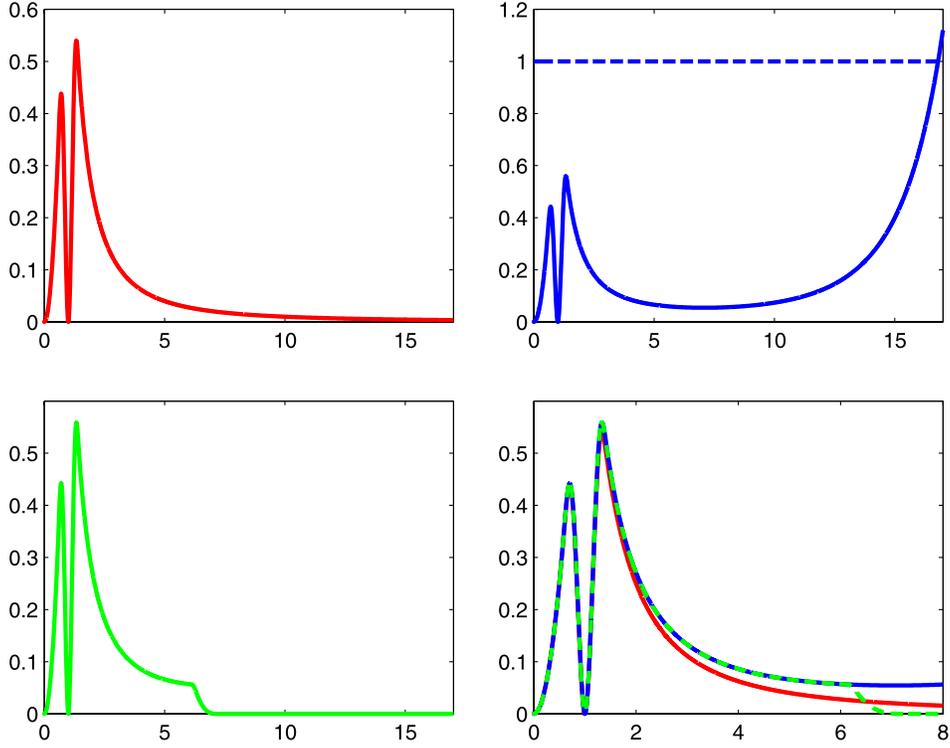

FIG. 1. *The first three panels illustrate $\hat{w}_1(t)$ (red), $\tilde{w}(t)$ (blue) and $\hat{w}_2(t)$ (green). Note that $\tilde{w}$ is not a characteristic function as $\tilde{w}(t) > 1$ for large $|t|$, and that $\hat{w}_2$ is a truncated version of $\tilde{w}$. The last panel is the overlay and zoom in of the first three panels.*

and

$$\varliminf_{n\to\infty} n^{2\beta} \cdot (\log n)^{(\alpha+1)} \cdot R_n^u(\mathcal{F}_0(\alpha,\beta,\varepsilon_0,q,a,A;n)) \geq C.$$

Due to the calibrations we choose in (2.3) and (2.5), the optimal rate is expressed in terms of parameters $\alpha, \beta$. Such calibrations are mainly for the simplicity in the presentation: Theorem 2.1 (as well as Theorems 2.2, 3.1 and 3.2 below) can be extended to more general settings. Here is an example. Fix $\varepsilon_0 \in (0,1)$ and $\beta \in [0, 1/2)$, suppose we (a) modify the calibration of $\varepsilon_n$ in (2.5) into that $\eta_n \leq \varepsilon_n \leq \varepsilon_0$ with $\eta_n$ being a sequence satisfying $\eta_n \geq \varepsilon_0 n^{-\beta}$, and (b) change the parameter space from $\mathcal{F}_0$ to $\mathcal{F}_0' = \mathcal{F}_0'(\alpha, q, a, A, \eta_n; n)$, where

$\mathcal{F}_0'(\alpha, \beta, q, a, A, \eta_n; n)$

$\quad = \{f \in \mathcal{F}$ and satisfies (2.2), (2.3), and constraints on $\varepsilon_n$ above$\}$.

The following corollary can be proved similarly as that of Theorem 2.1.

COROLLARY 2.1. *Fix $\alpha > 2$, $\beta \in [0, 1/2)$, $\varepsilon_0 \in (0,1)$, $q > 0$, $a > 0$ and $A > \sqrt{a^2 + 1} M_q^{1/q}$, let $\varepsilon_n$ and $\mathcal{F}_0'$ be calibrated as above. There is a constant $C > 0$ which depends on $\alpha, \beta, \varepsilon_0, q, a$ and $A$ such that*

$$\varliminf_{n \to \infty} \eta_n^{-2} \cdot (\log n)^{(\alpha+2)} \cdot R_n^\sigma(\mathcal{F}_0'(\alpha, \beta, q, a, A, \eta_n; n)) \geq C$$

*and*

$$\varliminf_{n \to \infty} \eta_n^{-2} \cdot (\log n)^{(\alpha+1)} \cdot R_n^u(\mathcal{F}_0'(\alpha, \beta, q, a, A, \eta_n; n)) \geq C.$$

We remark that for the case $\beta > 0$, the condition $A > \sqrt{a^2 + 1} M_q^{1/q}$ can be relaxed to that of $A > a M_q^{1/q}$. The latter is the minimum requirement for otherwise $\mathcal{F}_0(\alpha, \beta, \varepsilon_0, q, a, A; n)$ is an empty set. Theorem 2.1 shows that the minimax risk for estimating $\sigma_0^2$ cannot converge to 0 faster than $O(n^{-2\beta} \cdot (\log n)^{-(\alpha+2)})$, and that for estimating $u_0$ cannot be faster than $O(n^{-2\beta} \cdot (\log n)^{-(\alpha+1)})$. In next section, we shall show that these rates can indeed be attained and thus establish the minimax rates of convergence.

2.2. *Rate optimal estimators for the null parameters.* In this section, we seek estimators of the null parameters whose risks converge at the same rates as those of the lower bounds. Once such estimators are constructed, then their risks give upper bounds for the minimax risks, and the estimators themselves are rate optimal.

Given that estimating the null parameters is a relatively new problem, there are only a small number of methods in the literature. One straightforward approach is the method of moments, and another approach, proposed by Efron (2004), is to use the half-width of the central peak of the histogram. However, these approaches are only consistent in the sparse case where the proportion $\varepsilon = \varepsilon_n$ tends to 0 as $n$ tends to $\infty$. See Jin and Cai (2007) for more discussion.

In our recent work [Jin and Cai (2007)], we demonstrated that the null component can be well isolated in high-frequency Fourier coefficients, and based on this observation, we introduced a Fourier approach for estimating the null parameters. In detail, for any $t$ and complex-valued differentiable function $\xi$, let $\text{Im}(\xi)$ be the imaginary part and $\bar{\xi}$ be the complex conjugate, we introduce two functionals as follows:

(2.20)
$$\sigma_0^2(t; \xi) = -\left( \frac{d/ds |\xi(s)|}{s |\xi(s)|} \right) \bigg|_{s=t},$$
$$u_0(t; \xi) = \left( \frac{1}{|\xi(s)|^2} \cdot \text{Im}(\bar{\xi}(s) \xi'(s)) \right) \bigg|_{s=t}.$$

Next, fix $\gamma \in (0, 1/2)$, let $\varphi_n(t)$ be the *empirical characteristic function*,

$$\varphi_n(t) = \frac{1}{n} \sum_{j=1}^{n} e^{itX_j} \tag{2.21}$$

and

$$\hat{t}_n(\gamma) = \min\{t : t > 0, |\varphi_n(t)| \leq n^{-\gamma}\}. \tag{2.22}$$

We define the estimators for $\sigma_0^2$ and $u_0$ as

$$\hat{\sigma}_0^2(\gamma) = \sigma_0^2(\hat{t}_n(\gamma); \varphi_n), \qquad \hat{u}_0(\gamma) = u_0(\hat{t}_n(\gamma); \varphi_n).$$

To illustrate the idea behind the construction of these estimators, we consider a simplified case where $f$ is a homoscedastic Gaussian location mixture:

$$f(x) = (1 - \varepsilon) \frac{1}{\sigma_0} \phi\left(\frac{x - u_0}{\sigma_0}\right)$$
$$+ \varepsilon \int \frac{1}{\sigma_0} \phi\left(\frac{x - u}{\sigma_0}\right) h(u) \, du, \qquad h: \text{a univariate density}.$$

First, the empirical characteristic function approximates the *underlying characteristic function* $\varphi(t) = \varphi(t; f) \equiv E[e^{itX_j}]$,

$$\varphi_n(t) \approx \varphi(t) = e^{-\sigma_0^2 t^2 / 2}[(1 - \varepsilon) e^{iu_0 t} + \varepsilon \hat{h}(t)].$$

Second, by the well-known Riemann–Lebesgue lemma, for large $t$, $\hat{h}(t) \approx 0$, so

$$\varphi(t) \approx (1 - \varepsilon) e^{-\sigma_0^2 t^2 / 2} e^{iu_0 t} \equiv \varphi_0(t).$$

Last, $\hat{t}_n(\gamma)$ approximates its nonstochastic counterpart $t_n(\gamma)$,

$$t_n(\gamma) = \min\{t : t > 0, |\varphi(t)| \leq n^{-\gamma}\}. \tag{2.23}$$

Putting these together, we have that, heuristically,

$$\hat{\sigma}_0^2(\gamma) \approx \sigma_0^2(t_n(\gamma), \varphi_0) \equiv \sigma_0^2, \qquad \hat{u}_0(\gamma) \approx u_0(t_n(\gamma), \varphi_0) \equiv u_0,$$

where "$\equiv$" follow from direct calculations. See more discussions in Jin and Cai (2007).

The above approach has been studied in Jin and Cai (2007), where it was shown to be uniformly consistent across a wide class of cases. However, whether any of these estimators attains the optimal rate of convergence remains an open question. The difficulty is two-fold. First, compared to the study on consistency as in Jin and Cai (2007), the study on the optimal rate of convergence needs a much more delicate analysis on several small probability events. Tighter bounds on such events are not necessary for showing the consistency, but they are indispensable for proving the optimal rate of convergence. Second, a major technical difficulty is that the frequency $\hat{t}_n(\gamma)$ is stochastic and is not independent of the samples $X_j$. The stochasticity and dependence pose challenges in evaluating the estimation risks, and are the culprits for the lengthy analysis.

In this paper, we develop new analytical tools to solve these problems. The new analysis provides better probability bounds on several nuisance events and better control on the stochastic fluctuation of $\hat{t}_n(\gamma)$, $\hat{\sigma}_0^2(\gamma)$ and $\hat{u}_0(\gamma)$. The analysis reveals that the estimators $\hat{\sigma}_0^2(\gamma)$ and $\hat{u}_0(\gamma)$ are in fact rate-optimal under minimum regularity conditions. This is the following theorem, which is proved in Section 6.

THEOREM 2.2. *Fix* $\gamma \in (0, 1/2)$, $\alpha > 2$, $\beta \in [0, 1/2)$, $\varepsilon_0 \in (0,1)$, $q \geq 4$, $a > 0$ *and* $A > \sqrt{a^2 + 1} M_q^{1/q}$. *There is a constant* $C > 0$ *which only depends on* $\gamma, \alpha, \beta, \varepsilon_0, q, a$ *and* $A$ *such that*

$$\sup_{\mathcal{F}_0(\alpha,\beta,\varepsilon_0,q,a,A;n)} E[\hat{\sigma}_0^2(\gamma) - \sigma_0^2]^2 \leq C(n^{-2\beta} \log^{-(\alpha+2)}(n) + \log(n) \cdot n^{2\gamma-1})$$

*and*

$$\sup_{\mathcal{F}_0(\alpha,\beta,\varepsilon_0,q,a,A;n)} E[\hat{u}_0(\gamma) - u_0]^2 \leq C(n^{-2\beta} \log^{-(\alpha+1)}(n) + \log^2(n) \cdot n^{2\gamma-1}).$$

Taking $\gamma < 1/2 - \beta$ in Theorem 2.2, it then follows from Theorems 2.1 and 2.2 that the minimax rate of convergence for estimating the null parameters $\sigma_0^2$ and $u_0$ are $n^{-2\beta} \log^{-(\alpha+2)}(n)$ and $n^{-2\beta} \log^{-(\alpha+1)}(n)$, respectively. Furthermore, the estimators $\hat{\sigma}_0^2(\gamma)$ and $\hat{\mu}_0(\gamma)$ with $\gamma < 1/2 - \beta$ are rate optimal. Different choices of $\gamma$ does not affect the convergence rate but may affect the constant. In Section 4, we investigate how to choose $\gamma$ in practice with simulated data. We find that in many situations, the mean square error is relatively insensitive to the choice of $\gamma$, provided that it falls in the range of $(0.15, 0.25)$.

We mention that the logarithmic term in the minimax risk bears some similarity with the conventional deconvolution problem. See Section 5 for further discussion.

**3. Estimating the proportion of nonnull effects.** We now turn to the minimax estimation of the proportion. First, we consider the case where the null parameters are known. We show that, with careful modifications, the approach proposed in our earlier work [Jin and Cai (2007) and Jin (2008)] attains the optimal rate of convergence. We then extend the optimality to the case where the null parameters $(u_0, \sigma_0^2)$ are unknown.

3.1. *Estimating the proportion when the null parameters are known.* When the null parameters $(u_0, \sigma_0^2)$ are known, we can always use them to renormalize the test statistics $X_j$. So without loss of generality, we assume $u_0 = 0$ and $\sigma_0 = 1$. As a result, the marginal density of $X_j$ obeys a simplified form,

$$(3.1) \qquad X_j \stackrel{\text{i.i.d.}}{\sim} (1-\varepsilon)\phi(x) + \varepsilon \int \phi\left(\frac{x-u}{\sigma}\right) dH(u,\sigma) \equiv f.$$

The problem of estimating the proportion has received much recent attention. See, for example, Storey (2002), Genovese and Wasserman (2004), Meinshausen and Rice (2006) [see also Efron et al. (2001) and Swanepoel (1999)]. A necessary condition for the consistency of several of these approaches is that the marginal density of the nonnull effects (i.e., $f^{\text{alt}}$) is pure, a notion introduced in Genovese and Wasserman (2004). Unfortunately, the purity condition is generally not satisfied in the current setting; see Jin (2008) for a detailed discussion.

In our recent work Jin and Cai (2007) and Jin (2008), we proposed a Fourier approach to estimating the proportion which is described as follows. Let $\omega(\xi)$ be a bounded, continuous, and symmetric density function supported in $(-1, 1)$. Define a so-called phase function

$$\psi_n(t;\omega) = \psi_n(t;\omega,X_1,X_2,\ldots,X_n) = \int \omega(\xi) e^{t^2\xi^2/2} \varphi_n(t\xi)\,d\xi,$$

where as before $\varphi_n(t) = \frac{1}{n}\sum_{j=1}^n e^{itX_j}$ is the empirical characteristic function. Fix $\gamma \in (0, 1/2)$ and let $t_n = t_n(\gamma)$ be as in (2.23), the estimator is defined as

$$(3.2) \qquad \hat{\varepsilon}_n(\gamma;\omega) = \hat{\varepsilon}_n(\gamma;\omega,X_1,X_2,\ldots,X_n) = 1 - \text{Re}(\psi_n(t_n(\gamma);\omega)),$$

where $\text{Re}(z)$ stands for the real part of $z$. In Jin and Cai (2007) and Jin (2008), three different choices of $\omega(\xi)$ are recommended, namely the uniform density, the triangle density and the smooth density that is proportional to $\exp(-\frac{1}{1-|\xi|^2}) \cdot 1_{\{|\xi|<1\}}$.

The advantage of the Fourier approach is that it is no longer tied to the purity condition and can be shown to be consistent for the proportion uniformly for all eligible $H(u,\sigma)$; see details in Jin and Cai (2007) and Jin (2008). However, unfortunately, it is not hard to show that these estimators are not rate optimal with any of these three $\omega$.

In this paper, we propose the following estimator:

$$\hat{\varepsilon}_n(\gamma) = \left(1 - \frac{1}{n}\sum_{j=1}^n e^{t^2/2}\cos(tX_j)\right)\bigg|_{t=\sqrt{2\gamma \log n}}$$

$$= 1 - n^{-(1-\gamma)} \sum_{j=1}^n \cos(\sqrt{2\gamma \log n}\,X_j).$$

In comparison, $\hat{\varepsilon}_n(\gamma)$ is a special case of $\hat{\varepsilon}_n(\gamma;\omega)$, where instead of being a density function as in (3.2), $\omega$ is a point mass concentrated at 1. We shall show that under mild conditions, the proposed estimator $\hat{\varepsilon}_n(\gamma)$ attains the optimal rate of convergence. In detail, fix $\alpha > 0$, $\beta \in [0, 1/2)$, $\varepsilon_0 \in (0,1)$, $q \geq 2$, and $A > \sqrt{2} M_q^{1/2}$, let $\eta_n = \varepsilon_0 n^{-\beta}$ be as before. Consider the following parameter space for the minimax theory on estimating the proportion:

$$(3.3) \quad \tilde{\mathcal{F}} = \tilde{\mathcal{F}}(\alpha, \beta, \varepsilon_0, q, A; n) = \left\{ f \in \mathcal{F} : \varepsilon \leq \eta_n, \int |x|^q f(x)\,dx \leq A^q \right\}.$$

The minimax risk for estimating the proportion when the null parameters are known is

$$(3.4) \quad R_n^{\varepsilon,a} = R_n^{\varepsilon,a}(\tilde{\mathcal{F}}(\alpha, \beta, \varepsilon_0, q, A; n)) = \inf_{\hat{\varepsilon}} \left\{ \sup_{\tilde{\mathcal{F}}(\alpha,\beta,\varepsilon_0,q,A;n)} E[\hat{\varepsilon} - \varepsilon]^2 \right\}.$$

We have the following theorem.

THEOREM 3.1. *Fix $\gamma \in (0, 1/2)$, $\alpha > 0$, $\beta \in [0, 1/2)$, $\varepsilon_0 \in (0,1)$, $q \geq 2$, $A > \sqrt{2} M_q^{1/q}$. There is a generic constant $C > 0$ which only depends on $\alpha$, $\beta$, $\varepsilon_0$, $q$, $A$ and $\gamma$ such that for sufficiently large $n$,*

$$R_n^{\varepsilon,a}(\tilde{\mathcal{F}}(\alpha, \beta, \varepsilon_0, q, A; n)) \geq C n^{-2\beta} \log^{-\alpha}(n)$$

*and*

$$\sup_{\tilde{\mathcal{F}}_0(\alpha,\beta,\varepsilon_0,q,A;n)} E[\hat{\varepsilon}(\gamma) - \varepsilon]^2 \leq C(n^{-2\beta} \log^{-\alpha}(n) + n^{2\gamma - 1}).$$

*In particular, if $\gamma < 1/2 - \beta$, then $\hat{\varepsilon}_n(\gamma)$ attains the optimal rate of convergence.*

The proof of Theorem 3.1 is similar (but significantly simpler) than Theorem 3.2 below, which deals with the case where the null parameters are unknown. For reasons of space, we provide the proof of Theorem 3.2 in Section 6 but omit that of Theorem 3.1.

3.2. *Estimating the proportion when the null parameters are unknown.* We now turn to the case where the null parameters are unknown. A natural approach is to first estimate the null parameters with $(\hat{\sigma}_0(\gamma), \hat{u}_0(\gamma))$ and then plug them into $\hat{\varepsilon}_n(\gamma)$ to obtain an estimate of the proportion. In other words, fix $\gamma \in (0, 1/2)$, the plug-in estimator is

$$(3.5) \quad \hat{\varepsilon}_n^*(\gamma) = 1 - \frac{1}{n} \sum_{j=1}^n e^{t^2/2} \cos\left( t \left[ \frac{X_j - \hat{u}_0(\gamma)}{\hat{\sigma}_0(\gamma)} \right] \right) \bigg|_{\{t = \sqrt{2\gamma \log n}\}}.$$

We consider the minimax risk over the parameter space $\mathcal{F}_0$. The minimax risk for estimating the proportion when the null parameters are unknown is then

$$(3.6) \quad R_n^{\varepsilon,b} = R_n^{\varepsilon,b}(\mathcal{F}_0(\alpha,\beta,\varepsilon_0,q,A,a;n)) = \inf_{\hat{\varepsilon}} \left\{ \sup_{\mathcal{F}_0(\alpha,\beta,\varepsilon_0,q,a,A;n)} E[\hat{\varepsilon} - \varepsilon]^2 \right\}.$$

The following theorem, proved in Section 6, shows that the plug-in estimator is rate optimal.

THEOREM 3.2. *Fix $\gamma \in (0, 1/2)$, $\alpha > 2$, $\beta \in [0, 1/2)$, $\varepsilon_0 \in (0,1)$, $q > 4 + 2\gamma$, $a > 0$ and $A > \sqrt{a^2+1} M_q^{1/q}$. There is a generic constant $C > 0$ which only depends on $\gamma$, $\alpha$, $\beta$, $\varepsilon_0$, $q$, $a$ and $A$ such that for sufficiently large $n$,*

$$R_n^{\varepsilon,b} \geq C n^{-2\beta} \log^{-\alpha}(n)$$

*and*

$$\sup_{\mathcal{F}_0(\alpha,\beta,\varepsilon_0,q,a,A;n)} E[\hat{\varepsilon}_n^*(\gamma) - \varepsilon]^2 \leq C(n^{-2\beta} \log^{-\alpha}(n) + \log^3(n) \cdot n^{2\gamma-1}).$$

*Especially if $\gamma < 1/2 - \beta$, then $\hat{\varepsilon}_n^*(\gamma)$ attains the optimal rate of convergence.*

Compare Theorem 3.2 with Theorem 3.1, we see that except for the small difference in the upper bound [one has the $\log^3(n)$ term and the other does not], the minimax rates of convergence are the same whether the null parameters are known or not. The $\log^3(n)$ is the price we pay for the extra variability in estimation when the null parameters are unknown. Therefore, the plug-in estimator $\hat{\varepsilon}_n^*(\gamma)$ given in (3.5) is rate-optimal under almost the same conditions as in the case where the null parameters are known.

**4. Simulation study.** The procedures for estimating the proportion and null parameters presented in Sections 2 and 3 are easy to implement. In this section, we investigate the numerical performance of the procedure with simulated data.

The numerical study has several goals. The first is to consider the effect of the tuning parameter $\gamma$ on mean squared error (MSE) of the estimators and to make a recommendation on the choice of $\gamma$. The second is to compare the performance of the estimators with different $n$. The third is to compare the procedure with those in the literature. Several different combinations of the proportion and the nonnull distributions are used for such comparisons. The fourth is to investigate the performance of the estimators when the assumptions on eligibility and independence do not hold. The last and the most important goal is to study the effect of the estimation accuracy over the subsequent multiple testing procedures. Along this line, we consider two specific multiple testing procedures in our numerical study. One is the adaptive $p$-value based procedure (AP) introduced in Benjamini and Hochberg (2000) which requires an estimation of the proportion $\varepsilon$. This is the original Benjamini–Hochberg step-up procedure with an adjusted FDR level accounting for the sparsity. Another is the AdaptZ procedure (AZ) proposed in Sun and Cai (2007). This procedure thresholds the ranked Lfdr statistic (1.3) and requires estimations of $\varepsilon$, $f$ and $f^{\mathrm{null}}$. The procedure is asymptotically optimal in the sense that it minimizes the false nondiscovery rate asymptotically when the estimators of $\varepsilon$, $f$ and $f^{\mathrm{null}}$ are consistent.

Unless specified otherwise, the simulation results given in this section are based on $n = 10{,}000$, 1000 replications and the following Gaussian mixture model:

$$(4.1) \qquad X_i \sim (1-\varepsilon)N(\mu_0, \sigma_0^2) + \frac{\varepsilon}{2}N(\mu_{1i}, \sigma^2) + \frac{\varepsilon}{2}N(\mu_{2i}, \sigma^2),$$

where $\mu_{1i}$ and $\mu_{2i}$ are drawn from some distributions that may change from one case to another. Below, we report the simulation results along with the five aforementioned directions.

First, we study the effect of the tuning parameter $\gamma$ on the performance of the estimators. To this end, we consider the following setting.

SETTING 1. We take $\mu_0 = 0$, $\sigma_0 = 1$, $\mu_{1i} \sim \mathrm{Uniform}(-0.9, -0.1)$, $\mu_{2i} \sim \mathrm{Uniform}(0.5, 1.5)$, $\varepsilon = 0.2$ and $\sigma = 1.2$.

Table 1 tabulates the MSE of the three estimators $\hat{\varepsilon}_n^*(\gamma)$, $\hat{u}_0(\gamma)$ and $\hat{\sigma}_0^2(\gamma)$. The results suggest that $\hat{\varepsilon}_n^*$ and $\hat{\sigma}_0^2$ perform well in terms of the MSE when $\gamma$ is in a neighborhood of 0.2, ranging from 0.14 to 0.26 (note that, however, the estimator $\hat{\mu}_0$ favors a smaller $\gamma$). Additional simulations show similar patterns. In light of this, we conclude that an overall good choice is $\gamma = 0.2$. We recommend this choice for practical use in general, and use it in the rest of simulation study in this paper.

Second, we investigate how the number of hypotheses $n$ affects the estimation accuracy. The setting we consider is the same as Setting 1, but with different $n$.

SETTING 2. We take $\mu_0 = 0$, $\sigma_0 = 1$, $\mu_{1i} \sim \mathrm{Uniform}(-0.9, -0.1)$, $\mu_{2i} \sim \mathrm{Uniform}(0.5, 1.5)$, $\varepsilon = 0.2$, $\sigma = 1.2$ and $n$ ranges from 2000 to 500,000.

TABLE 1
*MSE (in unit of $10^{-4}$) of the estimators $\hat{\varepsilon}_n^*(\gamma)$, $\hat{u}_0(\gamma)$ and $\hat{\sigma}_0^2(\gamma)$ for different $\gamma$*

| $\gamma$ | 0.08 | 0.11 | 0.14 | 0.17 | 0.20 | 0.23 | 0.26 | 0.29 | 0.32 | 0.35 | 0.38 |
|---|---|---|---|---|---|---|---|---|---|---|---|
| MSE($\hat{\varepsilon}_n^*$) | 15.1 | 11.8 | 8.58 | 5.90 | 4.14 | 3.81 | 6.33 | 16.5 | 46.1 | 91.6 | 142 |
| MSE($\hat{u}_0$) | 0.37 | 0.93 | 1.79 | 3.11 | 5.40 | 9.65 | 17.8 | 33.3 | 63.0 | 114 | 204 |
| MSE($\hat{\sigma}_0^2$) | 2.31 | 1.57 | 1.07 | 0.78 | 0.68 | 0.77 | 1.08 | 1.70 | 2.83 | 4.89 | 8.84 |

TABLE 2
*Comparison of MSE (in unit of $10^{-5}$) for different n under Setting 2. The tuning parameter $\gamma$ is set at $0.2$*

| $n$ | 2000 | 5000 | 10,000 | 15,000 | 20,000 | 50,000 | 100,000 | 500,000 |
|---|---|---|---|---|---|---|---|---|
| $\text{MSE}(\hat{\varepsilon}_n^*)$ | 306.6 | 102.6 | 43.9 | 26.1 | 17.7 | 4.6 | 1.7 | 0.2 |
| $\text{MSE}(\hat{\mu}_0)$ | 596.6 | 143.8 | 60.5 | 31.7 | 19.3 | 5.8 | 1.9 | 0.2 |
| $\text{MSE}(\hat{\sigma}_0^2)$ | 74.6 | 19.6 | 7.1 | 3.95 | 2.5 | 0.6 | 0.2 | 0.01 |

Table 2 summarizes the MSE of the estimators under Setting 2. The results show that the accuracy of the estimators improves quickly as $n$ increases.

We now move to our third goal and compare the proposed estimator for the proportion with those in the literature, namely Efron's estimator $\hat{\varepsilon}^{\text{E}}$ [Efron (2004)] and Storey's estimator $\hat{\varepsilon}^{\text{S}}$ [Storey (2002), Genovese and Wasserman (2004)], assuming the null distribution is known. To distinguish from $\hat{\varepsilon}_n(\gamma)$, we denote the special case of $\gamma = 0.2$ by

$$\hat{\varepsilon}_n^{\text{CJ}} = \hat{\varepsilon}_n(0.2)$$

and may drop the subscript $n$ for simplicity. We compare these three estimators with data generated with different proportion $\varepsilon$ (Setting 3a) and different heteroscedasticity parameter $\sigma$ (Setting 3b).

SETTING 3a. We take $\mu_0 = 0$, $\sigma_0 = 1$, $\mu_{1i} \sim \text{Uniform}(-0.9, -0.1)$, $\mu_{2i} \sim \text{Uniform}(0.5, 1.5)$ and $\sigma = 1.2$. The value of $\varepsilon$ varies from 0.03 to 0.30. The goal is to see how the performance of the three estimators depends on the sparsity.

SETTING 3b. We set $\mu_0 = 0$, $\sigma_0 = 1$, $\mu_{1i} \sim \text{Uniform}(-0.9, -0.1)$, $\mu_{2i} \sim \text{Uniform}(0.5, 1.5)$ and $\varepsilon = 0.2$. The value of $\sigma$ varies from 1.2 to 2.1. The goal is to study the effect of the nonnull distribution on the estimation accuracy of the proportion estimators.

Table 3 tabulates the MSEs of these three point estimators. It is clear that our estimator $\hat{\varepsilon}^{\text{CJ}}$ performs well uniformly in all the cases. In particular it is robust under the various settings of nonnull distribution and sparsity. Table 3 shows that the MSE of $\hat{\varepsilon}^{\text{CJ}}$ increases gradually from $5.7 \times 10^{-5}$ to $10.1 \times 10^{-5}$ as $\varepsilon$ increases from 0.03 to 0.30. In comparison, the other two estimators $\hat{\varepsilon}^{\text{S}}$ and $\hat{\varepsilon}^{\text{E}}$ perform well in the sparse case but poorly in the nonsparse case. The MSEs of $\hat{\varepsilon}^{\text{E}}$ and $\hat{\varepsilon}^{\text{S}}$ increase about 120 times and 80 times, respectively, and they can sometimes be more than 10 times (some times even 39 times) larger than the MSE of $\hat{\varepsilon}^{\text{CJ}}$.

TABLE 3
*Comparison of MSE (in unit of $10^{-5}$) of three-point estimators $\hat{\varepsilon}^{\mathrm{CJ}}$, $\hat{\varepsilon}^{\mathrm{E}}$ and $\hat{\varepsilon}^{\mathrm{S}}$*

| | Setting 3a | | | | | | | | | |
|---|---|---|---|---|---|---|---|---|---|---|
| $\varepsilon$ | 0.03 | 0.06 | 0.09 | 0.12 | 0.15 | 0.18 | 0.21 | 0.24 | 0.27 | 0.30 |
| MSE($\hat{\varepsilon}^{\mathrm{CJ}}$) | 5.7 | 7.7 | 9.0 | 9.9 | 9.3 | 10.3 | 10.0 | 11.2 | 11.5 | 10.1 |
| MSE($\hat{\varepsilon}^{\mathrm{E}}$) | 3.3 | 14.6 | 33.4 | 60.3 | 95.8 | 139 | 190 | 249 | 316 | 394 |
| MSE($\hat{\varepsilon}^{\mathrm{S}}$) | 2.4 | 8.9 | 19.5 | 32.9 | 49.9 | 72.8 | 99.7 | 130 | 163 | 195 |
| | Setting 3b | | | | | | | | | |
| $\sigma$ | 1.2 | 1.3 | 1.4 | 1.5 | 1.6 | 1.7 | 1.8 | 1.9 | 2.0 | 2.1 |
| MSE($\hat{\varepsilon}^{\mathrm{CJ}}$) | 67.3 | 53.7 | 41.8 | 31.7 | 24.0 | 17.6 | 13.2 | 9.4 | 7.0 | 4.8 |
| MSE($\hat{\varepsilon}^{\mathrm{E}}$) | 172 | 164 | 153 | 146 | 138 | 129 | 122 | 114 | 108 | 100 |
| MSE($\hat{\varepsilon}^{\mathrm{S}}$) | 89.0 | 81.6 | 72.2 | 67.7 | 61.9 | 55.4 | 50.3 | 46.7 | 43.5 | 41.0 |

Next, we consider the case where either the assumption on eligibility or the assumption on independence is violated. Consider the eligible assumption first. Denote by $\mathrm{DE}(\mu, \tau)$ the double exponential distribution with the density function $f(x; \mu, \tau) = \frac{1}{2\tau} e^{-|x-\mu|/\tau}$. We shall generate $X_i$ as

$$(4.2) \qquad X_i \sim (1-\varepsilon) N(\mu_0, \sigma_0^2) + \frac{\varepsilon}{2} \mathrm{DE}(\mu_{1i}, \tau) + \frac{\varepsilon}{2} \mathrm{DE}(\mu_{2i}, \tau).$$

Since the double exponential can be viewed as a scale Gaussian mixture [West (1987)], it is seen that the eligible condition does not hold. Two different settings are considered.

SETTING 4a. We take $\mu_0 = 0$, $\sigma_0 = 1$ and assume the null parameters $\mu_0$ and $\sigma_0$ are known. First generate $\mu_{1i}$ from $U(-0.9, -0.1)$ and $\mu_{2i}$ from $U(0.5, 1.5)$, then generate $X_i$ as in (4.2) with $\tau = 1.2$. The proportion $\varepsilon$ varies from 0.03 to 0.30.

SETTING 4b. We take $\mu_0 = 0$, $\sigma_0 = 1$ and assume the null parameters $\mu_0$ and $\sigma_0$ are unknown. First generate $\mu_{1i}$ from $U(-0.9, -0.1)$ and $\mu_{2i}$ from $U(0.5, 1.5)$, then generate $X_i$ as in (4.2) with $\varepsilon = 0.2$. The value of $\tau$ varies from 1.2 to 2.1.

Table 4 gives the MSEs in Settings 4a and 4b. In Setting 4a, Efron's method is often found to be divergent numerically and is thus excluded from comparison. For small $\varepsilon$, Storey's method and our method yield similar results and both perform well. For moderate to large $\varepsilon$, however, our method demonstrates great superiority. In Setting 4b, Efron's method is again found to be divergent, and Storey's method does not apply as it requires the information of the null parameters. We therefore exclude both of them from comparison. In both settings, despite that the eligible condition is violated, our method continues to perform well.

TABLE 4
$MSE$ (in unit of $10^{-4}$) in Settings 4a and 4b

| Setting 4a | | | | | | | | | | |
|---|---|---|---|---|---|---|---|---|---|---|
| $\varepsilon$ | **0.03** | **0.06** | **0.09** | **0.12** | **0.15** | **0.18** | **0.21** | **0.24** | **0.27** | **0.30** |
| $\mathrm{MSE}(\hat{\varepsilon}^{\mathrm{CJ}})$ | 8.17 | 7.28 | 6.35 | 5.65 | 4.92 | 4.20 | 3.78 | 3.02 | 2.51 | 2.01 |
| $\mathrm{MSE}(\hat{\varepsilon}^{\mathrm{S}})$ | 3.25 | 6.79 | 9.76 | 14.35 | 19.93 | 19.69 | 23.68 | 21.67 | 21.01 | 20.18 |
| Setting 4b | | | | | | | | | | |
| $\tau$ | **1.2** | **1.3** | **1.4** | **1.5** | **1.6** | **1.7** | **1.8** | **1.9** | **2.0** | **2.1** |
| $\mathrm{MSE}(\hat{\varepsilon}_n^*)$ | 11.9 | 10.7 | 9.7 | 8.7 | 7.9 | 7.1 | 6.5 | 5.8 | 5.3 | 4.8 |
| $\mathrm{MSE}(\hat{\mu}_0)$ | 0.16 | 0.18 | 0.19 | 0.18 | 0.19 | 0.19 | 0.20 | 0.22 | 0.23 | 0.23 |
| $\mathrm{MSE}(\hat{\sigma}_0^2)$ | 4.1 | 4.1 | 4.2 | 4.2 | 4.0 | 3.9 | 3.7 | 3.6 | 3.5 | 3.3 |

The unsatisfactory behavior of Efron's estimator and Storey's estimator can be explained as follows. It is known in the literature that a necessary condition for Efron's estimator or Storey's estimator to be consistent is that the alternative density has a thinner tail than that of the null density either to the left or to the right [this is the so-called purity condition; see, e.g., Genovese and Wasserman (2004), Jin and Cai (2006) and Jin (2008)]. In Settings 4a and 4b, due to the heavy tail of the double exponential density, the purity condition is violated. It can be shown that asymptotically the bias of either Efron's estimator or Storey's estimator has the same magnitude as that of the true proportion. This explains why Efron's method does not always converge, and Storey's method has a reasonable performance when the underlying proportion is small, but behaves increasingly unsatisfactory as the proportion gets larger. This also suggests that, when the alternative density has a heavy tail, relying on the tail area for inference (as that in Efron's/Storey's method) can lead to a large bias. A promising alternative is the proposed Fourier-based method.

We now consider a case where the assumption on independence is violated. To do so, let $L$ be an integer that ranges from 0 to 50 with an increment of 10. For each $L$, we generate $n + L$ samples $w_1, w_2, \ldots, w_{n+L}$ from $N(0, 1)$, then let $z_j = \frac{1}{\sqrt{L+1}} \sum_{\ell=j}^{j+L} w_\ell$. The samples $z_j$ generated in this way are blockwise dependent with a block size $L$ (note that $L = 0$ corresponds the independent case). The setting we consider is as follows, where the null parameters are assumed as unknown.

SETTING 4c. Fix $\varepsilon = 0.2$ and $\sigma = 1.2$. Generate $X_i = z_i$ for $i = 1, 2, \ldots, 8000$, $X_i = \mu_{i1} + \sigma z_i$ for $8001 \leq i \leq 9000$, and $X_i = \mu_{i2} + \sigma z_i$ for $9001 \leq i \leq 10{,}000$, where $\mu_{1i}$ from $U(-0.9, -0.1)$ and $\mu_{2i}$ from $U(0.5, 1.5)$.

TABLE 5
*MSE (in unit of $10^{-3}$) in Setting 4c*

| $L$ | 0 | 10 | 20 | 30 | 40 | 50 |
|---|---|---|---|---|---|---|
| $\text{MSE}(\hat{\varepsilon}^{\text{CJ}})$ | 8.8 | 10.3 | 16.6 | 25.2 | 34.7 | 43.2 |
| $\text{MSE}(\hat{\mu}_0^{\text{CJ}})$ | 10.4 | 37.5 | 63.8 | 94.4 | 131.7 | 150.0 |
| $\text{MSE}(\hat{\sigma}_0^{\text{CJ}})$ | 5.4 | 13.5 | 23.0 | 34.8 | 49.3 | 52.1 |
| $\text{MSE}(\hat{\varepsilon}^{\text{E}})$ | 34.3 | 34.1 | 33.5 | 33.2 | 33.2 | 32.3 |
| $\text{MSE}(\hat{\mu}_0^{\text{E}})$ | 1.2 | 2.8 | 4.0 | 5.4 | 7.0 | 8.8 |
| $\text{MSE}(\hat{\sigma}_0^{\text{E}})$ | 14.7 | 18.1 | 21.7 | 28.1 | 34.7 | 33.5 |

Table 5 summarizes the results. In terms of MSE, the estimation accuracy decreases as the range of dependence increases. However, the MSE are still relatively small, especially those correspond to proportion and the null variance parameter $\sigma_0^2$. In comparison to Efron's method, correlation has a relatively larger impact on our method. The performance of our estimation procedure is better than Efron's when the correlation is weak to moderate. However, Efron's method is better when the correlation is strong.

The insight lies in the effect of correlation over the bias and variance. For all these estimators, the bias contains mainly marginal effects so the correlation does not have much effect on it. The correlation, however, may have important effect on the variance [see Jin and Cai (2006) and Jin (2008)]. In comparison, despite that our methods have a smaller bias, it gives relative larger MSE because it has a larger variance and is relatively more vulnerable when the correlation is strong.

Finally, we investigate how the point estimators affect the results of subsequent multiple testing procedures. First, we use the adaptive $p$-value based procedure [Benjamini and Hochberg (2000)] to compare the effect of the three point estimators of the proportion in the subsequent multiple testing. To this end, we consider the following two settings (which are the same as Setting 3a and 3b, respectively, but we restate them to avoid confusion).

SETTING 5a. We take $\mu_0 = 0$, $\sigma_0 = 1$, $\mu_{1i} \sim \text{Uniform}(-0.9, -0.1)$, $\mu_{2i} \sim \text{Uniform}(0.5, 1.5)$ and $\sigma = 1.2$. The value of $\varepsilon$ varies from 0.03 to 0.30.

SETTING 5b. We set $\mu_0 = 0$, $\sigma_0 = 1$, $\mu_{1i} \sim \text{Uniform}(-0.9, -0.1)$, $\mu_{2i} \sim \text{Uniform}(0.5, 1.5)$ and $\varepsilon = 0.2$. The value of $\sigma$ varies from 1.2 to 2.1.

It is known that the original step-up procedure of Benjamini and Hochberg (1995) is conservative: it controls the FDR level at $(1 - \varepsilon)\alpha$ instead of the nominal level $\alpha$. To remedy this shortcoming, Benjamini and Hochberg ([2000](#)) proposed an adaptive BH procedure which applies the original step-up procedure at level $\alpha' = \alpha/(1 - \hat{\varepsilon})$ instead of $\alpha$, where $\hat{\varepsilon}$ is an estimate of $\varepsilon$. Clearly the true FDR level of the adaptive BH procedure depends on the estimation accuracy of $\hat{\varepsilon}$.

We now compare the actual FDR level of the adaptive BH procedure using $\hat{\varepsilon}^{\text{CJ}}$, $\hat{\varepsilon}^{\text{S}}$, and $\hat{\varepsilon}^{\text{E}}$. In addition we also use the deviations of the false discovery proportion (FDP) from the nominal FDR level as a measure of the accuracy of the testing procedure. The FDP is a notion that is closely related to FDR: the FDP is the proportion of false positives among all rejections, and the FDR is the expected value of the FDP; see, for example, Genovese and Wasserman ([2004](#)). The deviations of the FDP from the nominal FDR level are naturally summarized by mean squared error. Denote the FDP of the adaptive BH procedure with the proportion being estimated by $\hat{\varepsilon}^{\text{E}}$, $\hat{\varepsilon}^{\text{S}}$ and $\hat{\varepsilon}^{\text{CJ}}$ by $\text{FDP}^{\text{E}}$, $\text{FDP}^{\text{S}}$ and $\text{FDP}^{\text{CJ}}$.

Figure [2](#) compares the actual FDR levels as well as the MSEs of $\text{FDP}^{\text{E}}$, $\text{FDP}^{\text{S}}$ and $\text{FDP}^{\text{CJ}}$. The two right panels are the ratios of the MSEs of $\text{FDP}^{\text{E}}$, $\text{FDP}^{\text{S}}$ and $\text{FDP}^{\text{CJ}}$ to $\text{MSE}(\text{FDP}^{\text{CJ}})$. In each of these settings, overall, the true FDR level of the adaptive BH procedure using $\hat{\varepsilon}^{\text{CJ}}$ is closest to the nominal level. The other two estimators, $\hat{\varepsilon}^{\text{E}}$ and $\hat{\varepsilon}^{\text{S}}$, tend to under-estimate the proportion $\varepsilon$ and consequently yield conservative testing procedure with the true FDR level below the nominal value. The FDP plots indicate that overall $\text{FDP}^{\text{E}}$ has larger deviations from the nominal FDR level in individual realizations than that of $\text{FDP}^{\text{S}}$ which is itself larger than that of $\text{FDP}^{\text{CJ}}$. These results show that our estimator $\hat{\varepsilon}^{\text{CJ}}$ yields the most accurate testing procedure: compared to $\text{FDP}^{\text{S}}$ and $\text{FDP}^{\text{E}}$, $\text{FDP}^{\text{CJ}}$ is not only smaller in biases, but also smaller in variances.

Next, we compare again our estimator of the null parameters with that by Efron ([2004](#)). But this time we do so by investigating the effect of different point estimators over the subsequent multiple testing procedures, namely the adaptiveZ procedure by Sun and Cai ([2007](#)). In detail, we consider the following setting.

SETTING 5c. We take $\mu_0 = 0$, $\mu_{1i} \sim \text{Uniform}(-0.9, -0.1)$, $\mu_{2i} \sim \text{Uniform}(0.5, 1.5)$, $\varepsilon = 0.2$, and $\sigma = 1.3$. The value of $\sigma_0$ varies from 0.5 to 1. In this setting we estimate both the proportion $\varepsilon$ and the null parameters $\mu_0$ and $\sigma_0$.

We now compare the performance of our estimators of the proportion and the null parameters with those of Efron ([2004](#)). [Storey ([2002](#)) assumed a known null distribution and did not provided estimators for the null parameters, so we exclude it from the comparison.] We compare the performance of these estimators as measured by the accuracy of the actual FDR level of the adaptive testing procedure introduced in Sun and Cai ([2007](#)). The AdaptZ procedure given in Sun and Cai (2007) aims to minimize the false nondiscovery rate subject to the constraint that the FDR level is controlled at a pre-specified level. This procedure thresholds the ordered Lfdr statistic

$$\widehat{\mathrm{Lfdr}}(z_i) = (1-\hat{\varepsilon})\tilde{f}^{\mathrm{null}}(z_i)/\tilde{f}(z_i),$$

where $\tilde{f}^{\mathrm{null}}$ and $\tilde{f}$ are estimators of $f^{\mathrm{null}}$ and $f$, respectively. The marginal density $f$ is estimated by a kernel density estimator with bandwidth chosen by cross-validation.

Figure 3 plots the true FDR levels of the AdaptZ procedure using our estimators of $\varepsilon$ and $\tilde{f}^{\mathrm{null}}$ with those of the same procedure using the estimators of $\varepsilon$ and $\tilde{f}^{\mathrm{null}}$ given in Efron (2004). Figure 3 also displays the ratio of the MSEs of the FDP of the two testing procedures, $\mathrm{MSE}(\mathrm{FDP}^{\mathrm{E}})/\mathrm{MSE}(\mathrm{FDP}^{\mathrm{CJ}})$. The results clearly show that the true FDR level of the testing procedure with our estimator is much closer to the nominal level than that with the estimators given in Efron (2004) and the FDP has smaller deviations from the nominal FDR level. Indeed, the $\mathrm{MSE}(\mathrm{FDP}^{\mathrm{CJ}})$ can sometimes be 15 times smaller than $\mathrm{MSE}(\mathrm{FDP}^{\mathrm{E}})$ [see Panel (b)].

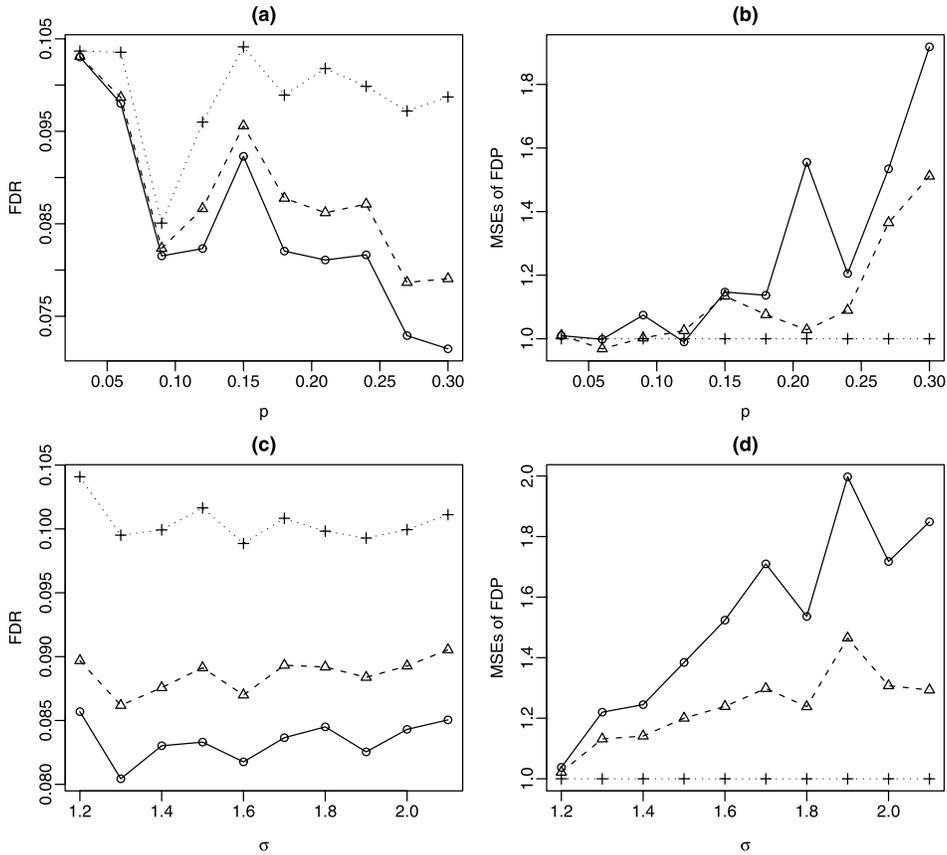

FIG. 2. *The actual FDR levels (left panels) and the MSEs of the FDP (right panels) of the adaptive BH procedure using the proportion estimators $\hat{\varepsilon}^{\mathrm{E}}$ (○ line), $\hat{\varepsilon}^{\mathrm{S}}$ (△ line) and $\hat{\varepsilon}^{\mathrm{CJ}}$ (+ line). The nominal level is 0.10. Top row: Setting 5a. The horizontal axis is the proportion $\varepsilon$. Bottom row: Setting 5b. The horizontal axis is the parameter $\sigma$.*

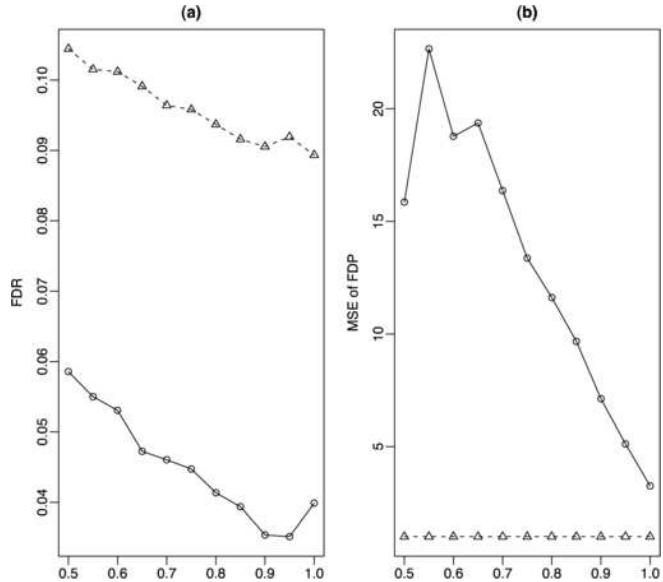

Fig. 3. *The actual FDR levels (left panel) and the relative MSEs of the FDP (right panel) of the AdaptZ procedure using the estimated null parameters and proportion: Efron's estimators (○ line) and our estimators (△ line). The nominal FDR level is* 0.10 *and the horizontal axis is the parameter* $\sigma_0$.

We conclude this section by mentioning that the proposed estimators usually yield a more accurate point estimation for the proportion and the null parameters than those by Efron (2004) and Storey (2002), not only asymptotically, but also for finite $n$. The accuracy of the proportion and the null parameters directly affects the performance of the subsequent testing procedures. Our estimators yield more accurate testing results than those by in Efron (2004) and Storey (2002).

**5. Discussion.** We derived the optimal rates of convergence for estimating the null parameters and the proportion of nonnull effects in large-scale multiple testing using a Gaussian mixture model. It was shown that the convergence rates depend on the smoothness of the mixing density $h(u|\sigma)$. The empirical characteristic function and Fourier analysis are crucial tools in our analysis of the optimality results. The proposed estimators not only are asymptotically rate-optimal but also enjoy superior finite $n$ performance.

Both theoretical and numerical results show that these estimators outperform the commonly used estimators in the literature. The improvement in the parameter estimation leads directly to more precise results in the subsequent multiple testing.

The minimax rates of convergence are proportional to the square of the true proportion multiplied by some logarithmic factors. The slowly convergent logarithmic factors can be attributed to the super-smooth nature of the Gaussian density, which attributes to the thin-tailed behavior of the corresponding characteristic function. As a result, even a relatively large perturbation in the true null parameters or in the true proportion may only result in a small difference in the $L^2$-norm of the characteristic function, which makes the perturbation hard to detect. The logarithmic terms are reminiscent of that found in the study of the conventional nonparametric deconvolution with Gaussian errors [e.g., Zhang (1990) and Fan (1991)], where the culprit for the slow convergence is also the super-smoothness of the Gaussian density. However, we should note that the problem considered here is different from the deconvolution problem; this explains the difference in the rate of minimax risk, the need for new procedures and the need for new approaches to derive the minimax risk bounds.

The work presented in this paper can be extended in several directions. First, while we have focused on the case where the characteristic function of $h$ decays at a polynomial rate, the results can be conveniently extended to the case where it has an exponential tail. Consider, for example, the following case:

$$|\hat{h}(t)| \leq C \exp(-|t|^\alpha).$$

The bias of the proposed estimator for the null parameters (and that for the proportion is similar) is of the order of

$$\exp(-C \log^{\alpha/2}(n)).$$

When $0 < \alpha < 2$, the bias is still larger than the variance and the rate of convergence is basically $\exp(-C \log^{\alpha/2}(n))$. When $\alpha > 2$, the bias tends to 0 faster than $1/\sqrt{n}$. In this case, the variance dominates the MSE, and we have $O(1/n)$ convergence rate. Second, while we focus on the case where $X_j$ are independent, extensions to the case of weak dependence is possible. Jin and Cai (2007) considered two dependent structures: the strongly $\alpha$-mixing case and the short-range dependent case and showed that the estimators constructed in that paper continue to be uniformly consistent under these dependent settings; see details therein. We expect that some of the results given in this paper are also extendable to the weakly dependent case. Third, while we focus on Gaussian mixtures in this paper, extensions to non-Gaussian mixtures is possible; see Jin (2008) for more discussion. An interesting example along this line is to replace the Gaussian mixture by the Laplace mixture. Due to the singularity of the Laplace density around the origin, the associated characteristic function decays much slower than that of the Gaussian density. As a result, the minimax risks for estimating the null parameters and the proportion are expected to have faster rates of convergence than those presented here. Last, while we focus on squared error loss here, the results can be extended to general loss functions.

We conclude this section by mentioning some possible future research directions. First, two key assumptions we make in this paper are the Gaussian mixture structure of the marginal density of the $z$-scores, and the independence among different $z$-scores. An interesting direction is to study the extend to which the presented results continue to hold when these assumptions are violated. An equally interesting direction is to study how to normalize/pre-process the data such that the assumptions hold approximately. Given the considerable efforts on normalization and pre-processing by the gene microarray community in recent years, the research along this direction could be very fruitful. Second, it would also be interesting to develop an adaptive approach to select the tuning parameter $\gamma$ in our proposed procedure. Given the overwhelming practical interest in large-scale multiple testing, this is an interesting problem for further study.

**6. Proof of the main results.** In this section, we prove the main results: Theorems 2.1, 2.2 and 3.2.

6.1. *Proof of Theorem 2.1.* The proofs of the minimax lower bounds for estimating the null parameters $\sigma_0^2$ and $u_0$ are similar. We present a detailed proof for the first claim and only a brief outline for the second one.

Consider the first claim. The key is to flesh out the ideas sketched in Section 2.1. We begin by filling in the details of the construction of $w_1$ and $w_2$. Let $k$ be the smallest even number that is greater than $2q+1$, let

$$\xi(t) = \begin{cases} \dfrac{(-1)^{k/2}\pi}{(k-1)!}|t|^{k-1}, & 0 \leq t \leq 1, \\ |t|^{-\alpha}, & t > 1, \end{cases}$$

and let $s_1$ and $s_2$ be two symmetric smooth functions, where $s_1$ satisfies (1). $0 \leq s_1(t) \leq 1$, (2). $s_1(t) = 1$ when $|t-1| > 2/3$, and (3). $s_1(t) = 0$ when $|t-1| < 1/3$, and $s_2$ satisfies (1). $0 \leq s_2(t) \leq 1$, (2). $s_2(t) = 1$ when $0 < |t| < \tau_n + 1/3$, and (3). $s_2(t) = 0$ when $|t| > \tau_n + 2/3$. The existence of such smooth function is well known in the literature; see Erdelyi (1956), for example. We construct $w_1$ and $w_2$ through their characteristic functions by

(6.1) $$\hat{w}_1(t) = s_1(t)\xi(t)$$

and

$$\hat{w}_2(t) = s_2(t) \cdot \left( e^{\delta_n t^2/2} \hat{w}_1(t) + \left( \frac{1}{\vartheta_0} \frac{1 - \eta_n}{\eta_n} \right) [e^{\delta_n t^2/2} - 1] \right); \tag{6.2}$$

see Figure 1 for illustrations.

Now, to show the claim, it remains to show (a) $h_1$ and $h_2$ are indeed densities, (b) the $\chi^2$-distance between $f_1$ and $f_2$ is equal to $o(1/n)$ and (c) the densities $f_1$ and $f_2$ in (2.9) and (2.10) satisfy the constraints (2.2) and (2.3) and therefore live in $\mathcal{F}_0(\alpha, \beta, \varepsilon_0, q, a, A; n)$. To do so, we need some lemmas.

Let $g$ be the Gaussian mixing density

$$g(x) = g(x; w_1, a) = \int \frac{1}{a} \phi\left(\frac{x-u}{a}\right) w_1(u).$$

By the way $f_1$ is defined in [see (2.9)], it is not hard to see that

$$f_1(x) = (1 - \eta_n)\phi_a(x) + \eta_n \phi_{\sqrt{a^2+1}}(x) + \vartheta_0 \eta_n g(x), \tag{6.3}$$

where $\phi_a$ denotes the density of $N(0, a^2)$. The following lemma characterizes the tail behavior of $w_1$, and so that of $g$ and $f_1$.

LEMMA 6.1. *For large $|u|$, $w_1(u) \sim |u|^{-k}$. As a result, for sufficiently small $\vartheta_0 > 0$, there is a constant $C > 0$,*

$$|g(x)| \leq C(1 + |x|)^{-k}, \qquad f_1(x) \geq C\eta_n(1 + |x|)^{-k}. \tag{6.4}$$

Here, $C > 0$ is a generic constant which only depends on (some or all) the parameters $\alpha, \beta, \varepsilon_0, q, a, A, k, \vartheta_0$ and $\theta_0$. The same rule applies below.

Next, the following lemma elaborates the tail behavior of $w_2$.

LEMMA 6.2. *For sufficiently large $|u|$ and $n$, there is a constant $C > 0$ such that*

$$||u|^k w_2(u) - 1| \leq C/|u|. \tag{6.5}$$

Last, the following lemma describes how close $f_1$ and $f_2$ are in the frequency domain.

LEMMA 6.3. *When $0 \leq |t| \leq \tau_n$, $\hat{f}_1(t) = \hat{f}_2(t)$. When $|t| > \tau_n$, there is a constant $C > 0$ such that for sufficiently large $n$,*

$$|\hat{f}_1^{(m)}(t) - \hat{f}_2^{(m)}(t)| \leq C|t|^m e^{-a^2 t^2/2}, \qquad m = 0, 1, \ldots, k/2.$$

Lemmas 6.1–6.3 are proved in the Appendix.

We are now ready to prove (a)–(c). Consider (a) first. By Lemmas 6.1 and 6.2, both $w_1$ and $w_2$ are positive for sufficiently large $|u|$. Therefore, (a) holds once we take $\vartheta_0$ sufficiently small.

Consider (b) next. Recall that the $\chi^2$-distance is $d(f_1, f_2) = \int [(f_1(x) - f_1(x))^2/f_1(x)] \, dx$. By (6.4) in Lemma 6.1,

$$\text{(6.6)} \quad \int \frac{|f_2(x) - f_1(x)|^2}{f_1(x)} \, dx \leq C\eta_n^{-1} \int (1 + |x|)^k |f_2(x) - f_1(x)|^2 \, dx$$
$$\leq Cn^\beta(\text{I} + \text{II}),$$

where $\text{I} = \int |f_2(x) - f_1(x)|^2 \, dx$ and $\text{II} = \int |x|^k (f_2(x) - f_1(x))^2 \, dx$. Now, by Parseval's formula [Mallat (1998)], for any integers $0 \leq m \leq k/2$,

$$\text{(6.7)} \quad \int x^{2m}|f_2(x) - f_1(x)|^2 \, dx = \int |x^m f_2(x) - x^m f_1(x)|^2 \, dx$$
$$= \int |\hat{f}_2^{(m)}(t) - \hat{f}_1^{(m)}(t)|^2 \, dt,$$

where by Lemma 6.3, the last term satisfies that

$$\text{(6.8)} \quad \int x^{2m}|f_2(x) - f_1(x)|^2 \, dx \leq C \int_{|t|>\tau_n} |t|^m e^{-a^2 t^2/2} \, dt.$$

Now, applying (6.8) to the case of $m = 0$ and $m = \ell$ gives

$$\text{(6.9)} \quad \text{I} + \text{II} \leq C \int_{|t|>\tau_n} (1 + |t|^{k/2}) e^{-a^2 t^2/2} \, dt \leq C\tau_n^{k/2-1} e^{-a^2 \tau_n^2/2}$$

and (b) follows by that $\beta < 1/2$ and that $a\tau_n = \sqrt{3 \log n}$.

Last, we show (c). It is sufficient to check both $f_1$ and $f_2$ satisfy (2.2) and (2.3). Consider $f_1$ first. Recall that $M_q$ is the $q$th moment of $N(0,1)$, combining (6.3) and (6.4) gives

$$\int |x|^q f(x) \, dx \leq [(1-\eta_n)a^q + \eta_n(a^2+1)^{q/2}]M_q + C\vartheta_0 \eta_n$$
$$\leq (a^2+1)^{q/2} M_q + C\vartheta_0 \varepsilon_0.$$

Therefore, by the assumption of $A > \sqrt{a^2+1} M_q^{1/q}$, (2.2) is satisfied once we take $\vartheta_0$ sufficiently small. At the same time, recall that $\hat{h}_1(t) = e^{-t^2/2} + \vartheta_0 \hat{w}_1(t)$ and that $\hat{w}_1(t) = |t|^{-\alpha}$ when $|t| \geq 4/3$, so (2.3) is also satisfied. Consider $f_2$ next. By Lemma 6.1 and the choice of $k$, the $2q$-moment of $f_1$ is finite. Using Hölder's inequality and (b),

$$\int |x|^q |f_1(x) - f_2(x)| \, dx \leq \left( \int |x|^{2q} f_1(x) \, dx \right)^{1/2} \left( \int \frac{(f_1(x) - f_2(x))^2}{f_1(x)} \, dx \right)^{1/2}$$
$$= o(1/\sqrt{n}).$$

Now, by the triangle inequality, $\int |x|^q f_2(x)\,dx \leq \int |x|^q f_1(x) + o(1/\sqrt{n})$, so $f_2$ satisfies the moment constraint in (2.2). At the same time, recall that $\hat{h}_2(t) = e^{-(1-\delta_n)t^2/2} + \vartheta_0 \hat{w}_2(t)$ and that

$$\hat{w}_2(t) = \begin{cases} e^{\delta_n t^2/2}\hat{w}_1(t) + \dfrac{1}{\vartheta_0}\dfrac{1-\eta_n}{\eta_n}[e^{\delta_n t^2/2} - 1], & |t| \leq \tau_n, \\ 0, & |t| \geq \tau_n + 1. \end{cases}$$

By elementary calculus and the choice of $\tau_n$ and $\delta_n$, there is a constant $C > 0$ such that for sufficiently large $n$ and $|t| > 4/3$,

$$|\hat{w}_2(t) - \hat{w}_1(t)| \leq C\theta_0 \tau_n^{-(\alpha+2)} t^2 \leq C\theta_0 |t|^{-\alpha},$$
$$|\hat{w}_2'(t) - \hat{w}_1'(t)| \leq C\theta_0 \tau_n^{-(\alpha+2)} t \leq C\theta_0 |t|^{-(\alpha+1)},$$

where we have used $w_1(t) = |t|^{-\alpha}$ for $|t| \geq 4/3$. Combining these we conclude that for a sufficiently small $\theta_0$, $h_2$ satisfies (2.3). This concludes the proof of (c) and the first claim of Theorem 2.1.

We now consider the second claim of Theorem 2.1. Similarly, the goal is to construct two density functions (say $f_3$ and $f_4$) in $\mathcal{F}_0(\alpha, \beta, \varepsilon_0, q, a, A; n)$ such that the null mean parameter $u_0$ associated with them differ by a small amount, and their $\chi^2$-distance is equal to $o(1/n)$. Let $\tau_n$, $s_2$, and $w_1$ be the same as in the proof associated with $\sigma_0^2$, and let $\theta_0 > 0$ be a constant to be determined. Define

(6.10)
$$\delta_n = \vartheta_0 \theta_0 \eta_n \tau_n^{-(\alpha+1)},$$
$$w_3 = w_1$$

and define $w_4$ through its characteristic function by

$$\hat{w}_4(t) = s_2(t) \cdot \left( \hat{w}_3(t) - \frac{2i}{\vartheta_0}\frac{1-\eta_n}{\eta_n}\sin\left(\frac{\delta_n t}{2}\right) \right).$$

We construct
$$\hat{h}_3(t) = e^{i\delta_n t/2}[e^{-t^2/2} + \vartheta_0 \cdot \hat{w}_3(t)], \qquad \hat{h}_4(t) = e^{-i\delta_n t/2}[e^{-t^2/2} + \vartheta_0 \cdot \hat{w}_4(t)],$$

and

(6.11) $\quad f_3(x) = (1-\eta_n)\dfrac{1}{a}\phi\left(\dfrac{x}{a}\right) + \eta_n \int \dfrac{1}{a}\phi\left(\dfrac{x-u}{a}\right) h_3(u)\,du,$

(6.12) $\quad f_4(x) = (1-\eta_n)\dfrac{1}{a}\phi\left(\dfrac{x-\delta_n}{a}\right) + \eta_n \int \dfrac{1}{a}\phi\left(\dfrac{x-\delta_n-u}{a}\right) h_4(u)\,du.$

Note that the null parameters associated with $f_3$ and $f_4$ differ by an amount of $\delta_n$. We are able to show that for appropriately small constants $\vartheta_0 > 0$ and $\theta_0 > 0$, $h_3$ and $h_4$ are indeed densities, and $f_3$ and $f_4$ live in $\mathcal{F}_0(\alpha, \beta, \varepsilon_0, q, a, A; n)$. Also, the $\chi^2$-distance between $f_3$ and $f_4$ is equal to $o(1/n)$. As the proofs are similar to that associated with $\sigma_0^2$, we skip them for reasons of space.

6.2. *Proof of Theorem 2.2.* Since the proofs are similar, we only prove the first claim. The following lemmas are proved in the Appendix.

LEMMA 6.4. *Fix $q \geq 4$ and $\gamma \in (0, 1/2)$. For sufficiently large $n$, and any event $B_n$ with $P\{B_n^c\} \leq C/n$, $E[\sigma_0^2(\varphi_n, \hat{t}_n(\gamma)) - \sigma_0^2)^2 \cdot 1_{\{B_n^c\}}] \leq C n^{2\gamma - 1}$.*

LEMMA 6.5. *Fix $q \geq 4$ and $\gamma \in (0, 1/2)$. For sufficiently large $n$,*
$$E[\varphi_n'(\hat{t}_n(\gamma)) - \varphi'(\hat{t}(\gamma))]^2 \leq C \log(n)/n.$$

We now proceed to show the theorem. Fix $q_1 > 3$, introduce the event

(6.13)
$$D_0 = \left\{ \frac{1}{n} \sum_{j=1}^n |X_j| \leq m_1 + 1, \right.$$
$$\left. \frac{1}{n} \sum_{j=1}^n X_j^2 \leq m_2 + 1, W_0(\varphi_n; n) \leq \sqrt{2q_1 \log n}/\sqrt{n} \right\},$$

where $m_1$ and $m_2$ are the first two moments of $X_1$ and
$$W_0(\varphi_n; n) = W_0(\varphi_n; n, X_1, X_2, \ldots, X_n) = \sup_{\{0 \leq t \leq \log n\}} |\varphi_n(t) - \varphi(t)|.$$

Note that, first, by Chebyshev's inequality,
$$P\left\{ \frac{1}{n} \sum_{j=1}^n |X_j| > m_1 + 1 \right\} \leq C/n, \qquad P\left\{ \frac{1}{n} \sum_{j=1}^n X_j^2 > m_2 + 1 \right\} \leq C/n.$$

Second, by Lemma A.2 of Jin and Cai (2007),
$$P\{W_0(\varphi_n; n) > \sqrt{2q_1 \log n}/\sqrt{n}\} \lesssim 4 \log^2(n) n^{-q_1/3}.$$

Recall that $q_1 > 3$, it thus follows that $P\{D_n^c\} \leq C/n$. By Lemma 6.4, $D_n^c$ only has a negligible contribution to the mean squared errors:

(6.14) $$E[(\sigma_0^2(\varphi_n, \hat{t}_n(\gamma)) - \sigma_0^2)^2 \cdot 1_{\{D_0^c\}}] \leq C n^{2\gamma - 1}$$

and all remains to show is

(6.15) $$E[(\sigma_0^2(\varphi_n, \hat{t}_n(\gamma)) - \sigma_0^2)^2 \cdot 1_{\{D_0\}}]$$
$$\leq C[n^{-2\beta} \log^{-(\alpha+2)}(n) + \log(n) n^{2\gamma - 1}].$$

We now show (6.15). Write for short $\hat{t}_n = \hat{t}_n(\gamma)$ and $t_n = t_n(\gamma)$. By the triangle inequality,
$$|\sigma_0^2(\varphi_n, \hat{t}_n) - \sigma_0^2| \leq |\sigma_0^2(\varphi_n, \hat{t}_n) - \sigma_0^2(\varphi, \hat{t}_n)| + |\sigma_0^2(\varphi, \hat{t}_n) - \sigma_0^2(\varphi, t_n)|$$
$$+ |\sigma_0^2(\varphi, t_n) - \sigma_0^2|.$$

So to show (6.15), all we need to show are

$$(6.16) \qquad E[(\sigma_0^2(\varphi_n, \hat{t}_n) - \sigma_0^2(\varphi, \hat{t}_n))^2 \cdot 1_{\{D_0\}}] \leq C \log(n) n^{2\gamma-1},$$

$$(6.17) \qquad E[(\sigma_0^2(\varphi, \hat{t}_n) - \sigma_0^2(\varphi, t_n))^2 \cdot 1_{\{D_0\}}] \leq C n^{2\gamma-1}$$

and

$$(6.18) \qquad |\sigma_0^2(\varphi, t_n) - \sigma_0^2| \leq C n^{-\beta} \log^{-(\alpha+2)/2}(n) \qquad \text{over } D_0.$$

Below, we show (6.16)–(6.18) separately.

Consider (6.16) first. By Lemmas A.2 and A.3 of Jin and Cai (2007), over the event $D_0$,

$$(6.19) \quad |\varphi_n(\hat{t}_n) - \varphi(\hat{t}_n)| \leq C\sqrt{\log n}/\sqrt{n}, \qquad |\hat{t}_n - t_n| \leq c_0 n^{\gamma-1/2},$$

where $c_0 > \sigma_0 \sqrt{q_1/\gamma}$ is a constant. Apply Lemma 6.1 of Jin and Cai (2006) with $f = \varphi_n$, $g = \varphi$, and $t = \hat{t}_n$,

$$(6.20) \quad \begin{aligned} &|\sigma_0^2(\varphi_n, \hat{t}_n) - \sigma_0^2(\varphi, \hat{t}_n)| \\ &\lesssim n^\gamma \bigg[ 3\sigma_0^2 |\varphi_n(\hat{t}_n) - \varphi(\hat{t}_n)| + \frac{1}{\hat{t}_n} |\varphi_n'(\hat{t}_n) - \varphi'(\hat{t}_n)| \bigg]. \end{aligned}$$

Combining (6.19) and (6.20) gives that, over the event $D_0$,

$$|\sigma_0^2(\varphi_n, \hat{t}_n) - \sigma_0^2(\varphi, \hat{t}_n)| \leq C n^\gamma \bigg[ \frac{\sqrt{\log n}}{\sqrt{n}} + \frac{1}{t_n} |\varphi_n'(\hat{t}_n) - \varphi'(\hat{t}_n)| \bigg]$$

and applying the Lemma 6.5 gives (6.16).

Consider (6.17) next. Direct calculations show that $|\frac{d}{dt}\sigma_0^2(\varphi, t)| \leq C$ for sufficiently large $t$. Using the second part of (6.19),

$$(6.21) \quad |\sigma_0^2(\varphi, \hat{t}_n) - \sigma_0^2(\varphi, t_n)| \leq C|\hat{t}_n - t_n| \leq C n^{\gamma-1/2} \qquad \text{over } D_0,$$

and (6.17) follows directly.

Last, we consider (6.18). Similar to Lemma 6.5 of Jin and Cai (2007), $|\sigma_0^2(\varphi, t_n) - \sigma_0^2| \leq C \frac{|\psi'(t_n)|}{t_n}$, where $\psi(t) = \varepsilon_n \int e^{it(u-u_0) - (\sigma^2-\sigma_0^2)t^2/2} \times h(u|\sigma) \, dH(\sigma)$. By direct calculations,

$$\begin{aligned} |\psi'(t)| &= \varepsilon_n \bigg| \int (i(u-u_0) - (\sigma^2-\sigma_0^2)t) e^{it(u-u_0)-(\sigma^2-\sigma_0^2)t^2/2} h(u|\sigma) \, dH(\sigma) \bigg| \\ &\leq \mathrm{I} + \mathrm{II}, \end{aligned}$$

where

$$\mathrm{I} = \varepsilon_n \bigg| \int (u-u_0) e^{it(u-u_0)-(\sigma^2-\sigma_0^2)t^2/2} h(u|\sigma) \, dH(\sigma) \bigg|$$

and

$$\mathrm{II} = \varepsilon_n \left| \int (\sigma^2 - \sigma_0^2) t e^{it(u-u_0) - (\sigma^2 - \sigma_0^2)t^2/2} h(u|\sigma) \, dH(\sigma) \right|.$$

By elementary Fourier analysis and the definition of $\hat{h}(t|\sigma)$ and $\tilde{h}(t|\sigma)$ [see (2.4)],

$$\mathrm{I} = \varepsilon_n \left| \int \tilde{h}'(t|\sigma) e^{-(\sigma^2 - \sigma_0^2)t^2/2} \, dH(\sigma) \right| \leq \varepsilon_n \int |\tilde{h}'(t|\sigma)| \, dH(\sigma)$$

and

$$\mathrm{II} \leq \varepsilon_n \left| \int \hat{h}(t|\sigma)(\sigma^2 - \sigma_0^2) t e^{-(\sigma^2 - \sigma_0^2)t^2/2} \, dH(\sigma) \right| \leq C(\varepsilon_n/t) \int |\hat{h}(t|\sigma)| \, dH(\sigma),$$

where we have used the fact that $\sup_{a \geq 0}\{ate^{-at^2/2}\} \leq C/t$. Combining these with (2.3) and (2.5) gives (6.18). This concludes the proof of Theorem 2.2.

6.3. *Proof of Theorem 3.2.* Consider the first claim first. Similar to the construction of the minimax lower bound for estimating the null parameter $\sigma_0^2$, the goal is to construct two density functions (say $f_5$ and $f_6$) in $\mathcal{F}_0(\alpha, \beta, \varepsilon_0, q, a, A; n)$ such that the proportion associated with them differ by a small amount, and their $\chi^2$-distance is equal to $o(1/n)$.

We construct $f_5$ and $f_6$ as follows. Let $\tau_n$, $w_1$, and $s_2$ be the same as in Section 6.1. Similarly, for a constant $\theta_0 > 0$ to be determined, let

$$\delta_n = \vartheta_0 \theta_0 \eta_n \tau_n^{-\alpha},$$
(6.22)
$$w_5 \equiv w_1$$

and

$$\hat{w}_6(t) = s_2(t) \cdot \left( \frac{\eta_n - \delta_n}{\eta_n} \hat{w}_5(t) + \frac{1}{\vartheta_0} \frac{\delta_n}{\eta_n}(1 - e^{-t^2/2}) \right).$$

We define $h_5$ and $h_6$ as

$$\hat{h}_5(t) = e^{-t^2/2} + \vartheta_0 \cdot \hat{w}_5(t), \qquad \hat{h}_6(t) = e^{-t^2/2} + \vartheta_0 \cdot \hat{w}_6(t),$$

and

(6.23) $$f_5(x) = (1 - \eta_n + \delta_n)\frac{1}{a}\phi\left(\frac{x}{a}\right) + (\eta_n - \delta_n)\int \frac{1}{a}\phi\left(\frac{x-u}{a}\right)h_5(u)\,du,$$

(6.24) $$f_6(x) = (1 - \eta_n)\frac{1}{a}\phi\left(\frac{x}{a}\right) + \eta_n \int \frac{1}{a}\phi\left(\frac{x-u}{a}\right)h_6(u)\,du.$$

Note that the proportion associated with $f_5$ and $f_6$ differ by an amount of $\delta_n$. We are able to show that for appropriately small constants $\vartheta_0 > 0$ and $\theta_0 > 0$, $h_5$ and $h_6$ are indeed densities, and $f_5$ and $f_6$ live in $\mathcal{F}_0(\alpha, \beta, \varepsilon_0, q, a, A; n)$. Also, the $\chi^2$-distance between $f_5$ and $f_6$ is equal to $o(1/n)$. As the proofs are similar to the case for $\sigma_0^2$, we skip them for reasons of space.

We now consider the second claim. Write for short $\hat{\varepsilon}_n^* = \hat{\varepsilon}_n^*(\gamma)$, $\hat{\sigma}_0^2 = \sigma_0^2(\gamma)$, and $\hat{u}_0 = \hat{u}_0(\gamma)$, introduce the nonstochastic counterparts of $\hat{\sigma}_0^2$ and $\hat{u}_0$, respectively, by

$$\bar{\sigma}_0^2 = \sigma_0(\varphi, t_n), \qquad \bar{u}_0 = u_0(\varphi, t_n),$$

where $t_n$ is defined in (2.23). The following lemma is a direct result of Theorem 1 of Jin and Cai (2007), which elaborates the stochastic fluctuation of $\hat{\sigma}_0^2$ and $\hat{u}_0$.

LEMMA 6.6. *Let $\gamma \in (0, 1/2)$ and $q > 4 + 2\gamma$ be as in the theorem. There is an event $B_n$ such that $P\{B_n^c\} = o(1/n)$ and over the event $B_n$,*

$$(6.25) \quad |\hat{\sigma}_0^2 - \bar{\sigma}_0^2| \leq C \log^{1/2}(n) n^{\gamma - 1/2}, \qquad |\hat{u}_0 - \bar{u}_0| \leq C \log(n) n^{\gamma - 1/2}.$$

Now, by replacing $\hat{u}_0$ with $\bar{u}_0$ in the definition of $\hat{\varepsilon}_n^*$, we introduce the following pseudo-estimator:

$$(6.26) \quad \tilde{\varepsilon}_n = \tilde{\varepsilon}_n(\gamma, X_1, \ldots, X_n, \bar{u}_0) = 1 - n^{\gamma - 1} \sum_{j=1}^n \cos\left(\sqrt{2\gamma \log n} \frac{X_j - \bar{u}_0}{\hat{\sigma}_0}\right).$$

The pseudo-estimator plays a key role in the proof. To see the point, we need some notation. Let $\tilde{\varphi}_n$ be the empirical characteristic function corresponding to $(X_j - \bar{u}_0)/\bar{\sigma}_0$,

$$(6.27) \qquad \tilde{\varphi}_n(t) = \tilde{\varphi}_n(t; X_1, \ldots, X_n; \bar{u}_0, \bar{\sigma}_0) = \frac{1}{n} \sum_{j=1}^n e^{it(X_j - \bar{u}_0)/\bar{\sigma}_0},$$

let $\tilde{\varphi}(t)$ be the corresponding (underlying) characteristic function

$$\tilde{\varphi}(t) = \tilde{\varphi}(t; f, \bar{u}_0, \bar{\sigma}_0) \equiv E[\tilde{\varphi}_n(t)]$$

and denote the real part of $\tilde{\varphi}_n$ and $\tilde{\varphi}$ by $\tilde{\varphi}_n^R$ and $\tilde{\varphi}^R$, respectively. Observe that if we denote

$$(6.28) \qquad \tilde{t}_n = \tilde{t}_n(\gamma; \hat{\sigma}_0, \bar{\sigma}_0) = \frac{\bar{\sigma}_0}{\hat{\sigma}_0} \sqrt{2\gamma \log n},$$

then $\tilde{\varepsilon}_n$ can be rewritten as

$$(6.29) \qquad \tilde{\varepsilon}_n = 1 - n^\gamma \tilde{\varphi}_n^R(\tilde{t}_n).$$

The advantage of introducing $\tilde{\varepsilon}_n$ is two-fold. First, by elementary trigonometrics, the difference between $\hat{\varepsilon}_n^*$ and $\tilde{\varepsilon}$ has a very simple form. This is the following lemma, whose proof is elementary so we omit it.

LEMMA 6.7.
$$\hat{\varepsilon}_n^* - \tilde{\varepsilon}_n = n^\gamma \operatorname{Re}\left(\tilde{\varphi}_n(\tilde{t}_n) \cdot \left[\sin^2\left(\frac{\tilde{t}_n}{2}\frac{\bar{u}_0 - \hat{u}_0}{\bar{\sigma}_0}\right) - i\sin\left(\tilde{t}_n\frac{\bar{u}_0 - \hat{u}_0}{\bar{\sigma}_0}\right)\right]\right).$$

Second, the stochastic fluctuation of $\tilde{\varepsilon}_n$ can be conveniently bounded through the maximum deviation of $\tilde{\varphi}_n(t)$ over the interval, say, $[0, \log(n)]$. In detail, fix a constant $q_1 > 3$, introduce the following event:
$$\tilde{D}_0 = \left\{\sup_{\{0 \leq t \leq \log n\}}\{|\tilde{\varphi}_n(t) - \tilde{\varphi}(t)|\} \leq \sqrt{2q_1 \log n}/\sqrt{n}\right\}.$$

The following lemma can be proved similarly as that of Lemma A.2 in Jin and Cai (2007), so we omit it.

LEMMA 6.8. $P\{\tilde{D}_0^c\} \lesssim 4\log^2(n)n^{-q_1/3}$.

A direct consequence of Lemma 6.8 is that
$$E|\tilde{\varphi}_n(\tilde{t}_n) - \tilde{\varphi}(\tilde{t}_n)|^2 \leq E\left[\sup_{\{0 \leq t \leq \log n\}}\{|\tilde{\varphi}_n(t) - \tilde{\varphi}(t)|^2\}\right] + o(1/n)$$
(6.30)
$$\leq C\log(n)/n.$$

Given the lemmas above, what remains to analyze is $\tilde{\varphi}^R(\tilde{t}_n)$. Note that $\tilde{t}_n$ fluctuates around $\sqrt{2\gamma \log n}$. We have the following lemma, which is proved in the Appendix.

LEMMA 6.9. Let $B_n$ be the event as in Lemma 6.6. We have
$$|\tilde{\varphi}^R(\tilde{t}_n) - \tilde{\varphi}^R(\sqrt{2\gamma \log n})| \leq C\log^{3/2}(n)n^{-1/2} \qquad \text{over } B_n$$
and
$$|(1 - \varepsilon_n) - n^\gamma \tilde{\varphi}^R(\sqrt{2\gamma \log n})| \leq C\log^{-\alpha/2}(n)n^{-\beta}.$$

We are now ready to show the theorem. By the triangle inequality and the Cauchy–Schwarz inequality,
$$|\hat{\varepsilon}_n^* - \varepsilon_n|^2$$
(6.31)
$$\leq (|\hat{\varepsilon}_n^* - \tilde{\varepsilon}_n| + |\tilde{\varepsilon}_n - (1 - n^\gamma \tilde{\varphi}^R(\tilde{t}_n))| + |(1 - n^\gamma \tilde{\varphi}^R(\tilde{t}_n)) - \varepsilon_n|)^2$$
$$\leq C(|\hat{\varepsilon}_n^* - \tilde{\varepsilon}_n|^2 + |\tilde{\varepsilon}_n - (1 - \tilde{\varphi}^R(\tilde{t}_n))|^2 + |(1 - n^\gamma \tilde{\varphi}^R(\tilde{t}_n)) - \varepsilon_n|^2).$$

First, by (6.29) and (6.30),

(6.32) $E|\tilde{\varepsilon}_n - (1 - n^\gamma \tilde{\varphi}^R(\tilde{t}_n))|^2 = n^{2\gamma}E|\tilde{\varphi}_n^R(\tilde{t}_n) - \tilde{\varphi}^R(\tilde{t}_n)|^2 \leq C\log(n)n^{2\gamma-1}.$

Second, by Lemma 6.9 and the Cauchy–Schwarz inequality,

$$
(6.33) \quad \begin{aligned}
E|(1-\varepsilon_n) - n^\gamma \tilde{\varphi}^R(\tilde{t}_n)|^2 &\leq C[\log^{-\alpha/2}(n)n^{-\beta} + \log^{3/2}(n)n^{\gamma-1/2}]^2 \\
&\leq C[\log^{-\alpha}(n)n^{-2\beta} + \log^3(n)n^{2\gamma-1}].
\end{aligned}
$$

Plugging this into (6.31) gives

$$(6.34) \quad E|\hat{\varepsilon}_n^* - \varepsilon_n|^2 \leq C[E|\hat{\varepsilon}_n^* - \tilde{\varepsilon}_n|^2 + \log^{-\alpha}(n)n^{-2\beta} + \log^3(n)n^{2\gamma-1}].$$

Compare (6.34) with the theorem. All that remains to show is

$$(6.35) \quad E|\hat{\varepsilon}_n^* - \tilde{\varepsilon}_n|^2 \leq C\log^3(n)n^{2\gamma-1}.$$

We now show (6.35). Note that $|\hat{\varepsilon}_n^* - \tilde{\varepsilon}_n|^2 \leq n^{2\gamma}$ and $P\{\tilde{D}_0^c \cup B_n^c\} = o(1/n)$, so

$$E[|\hat{\varepsilon}_n^* - \tilde{\varepsilon}_n|^2 \cdot 1_{\{\tilde{D}_0^c \cup B_n^c\}}] \leq o(n^{2\gamma-1})$$

and all we need to show is

$$(6.36) \quad E[|\hat{\varepsilon}_n^* - \tilde{\varepsilon}_n|^2 \cdot 1_{\{\tilde{D}_0 \cap B_n\}}] \leq C\log^3(n)n^{2\gamma-1}.$$

To this end, note that over the event $\tilde{D}_0 \cap B_n$, by Lemma 6.7 and that $|\sin(x)| \leq C|x|$ for all $x$,

$$(6.37) \quad |\hat{\varepsilon}_n^* - \tilde{\varepsilon}_n|^2 \leq C\tilde{t}_n^2 |\varphi_n(\tilde{t}_n)|^2 \frac{(\hat{u}_0 - \bar{u}_0)^2}{\bar{\sigma}^2}.$$

Now, first, by Lemma 6.6,

$$(6.38) \quad \tilde{t}_n \sim \sqrt{2\log n}, \qquad |\hat{u}_0 - \bar{u}_0| \leq C\log(n)n^{\gamma-1/2}(n).$$

Second, by Lemma 6.8 and the Cauchy–Schwarz inequality,

$$|\varphi_n(\tilde{t}_n)|^2 \leq \left|\tilde{\varphi}(\tilde{t}_n) + \frac{\sqrt{2q_1 \log n}}{\sqrt{n}}\right|^2,$$

where according to Lemma 6.9,

$$\varphi(\tilde{t}_n) \leq Cn^{-\gamma}.$$

Therefore, over the event $\tilde{D}_0 \cap B_n$,

$$(6.39) \quad |\varphi_n(\tilde{t}_n)|^2 \leq Cn^{-2\gamma}.$$

Inserting (6.38) and (6.39) into (6.37) gives (6.36), and concludes the proof of the theorem.

## APPENDIX

We shall prove in this section the technical lemmas which are used in the proofs of the main results in the previous sections.

**A.1. Proof of Lemma 6.1.** Consider the first claim first. The symmetry of $\hat{w}$ implies
$$w_1(u) = \frac{1}{2\pi}\int e^{-itu}\hat{w}_1(t)\,dt = \frac{1}{\pi}\int_0^\infty \cos(tu)\hat{w}_1(t)\,dt.$$
Note that $\hat{w}_1$ is smooth in $(0,\infty)$ and $\hat{w}_1^{(k-1)}(0) = (-1)^{k/2}\pi$. Repeatedly using integration by parts $k$ times yields
$$(A.1) \qquad \frac{1}{\pi}\int_0^\infty \cos(tu)\hat{w}_1(t)\,dt = u^{-k} + r_1(u), \qquad u > 0,$$
where
$$|r_1(u)| = \frac{1}{\pi|u|^k}\left|\int_0^\infty \cos(tu)\hat{w}_1^{(k)}(t)\,dt\right| = \frac{1}{\pi|u|^{k+1}}\left|\int_0^\infty \sin(tu)\hat{w}_1^{(k+1)}(t)\,dt\right|.$$
Direct calculations show that there is a constant $C = C(\alpha,k) > 0$ such that
$$|\hat{w}_1^{(k+1)}(t)| \leq C(1+|t|)^{-(\alpha+k+1)},$$
so
$$(A.2) \qquad |r_1(u)| \leq C|u|^{-(k+1)}.$$
Combining (A.1) and (A.2) gives the claim.

Next, consider the second claim. It is sufficient to show that for sufficiently large $x$,
$$(A.3) \qquad g(x) \geq C|x|^{-k}.$$
By the way $g$ is defined,
$$(A.4) \qquad g(x) = \int \phi_a(x)w_1(x-u)\,du = \mathrm{I} + \mathrm{II},$$
where
$$\mathrm{I} = \int_{|u|\geq x/2} w_1(x-u)\phi_a(u)\,du, \qquad \mathrm{II} = \int_{|u|<x/2} w_1(x-u)\phi_a(u)\,du.$$
First, we have
$$(A.5) \qquad |\mathrm{I}| \leq C\phi_a(x/2).$$
Second, by the first claim, there are generic constants $C_2 > C_1 > 0$ such that for sufficiently large $x$ and $|u| < x/2$,
$$C_1|x|^{-k} \leq w_1(x-u) \leq C_2|x|^{-k},$$
and so
$$(A.6) \qquad C_1(1+|x|)^{-k} \leq \mathrm{II} \leq C_2|x|^{-k}.$$
Inserting (A.5) and (A.6) into (A.4) gives (A.3).

Last, consider the third claim. Recall that [i.e., (6.3)]
$$f_1(x) = (1-\eta_n)\phi_a(x) + \eta_n\phi_{\sqrt{a^2+1}}(x) + \vartheta_0\eta_n g(x).$$
Once we take $\vartheta_0$ appropriately small, the claim follows from (A.3).

**A.2. Proof of Lemma 6.2.** Similarly, write

$$w_2(u) = \frac{1}{\pi} \int \cos(tu) \hat{w}_2(u) \, du = u^{-k} + r_2(u), \qquad u > 0,$$

where

$$|r_2(u)| \leq \frac{1}{\pi |u|^{k+1}} \int_0^\infty |\hat{w}_2^{(k+1)}(t)| \, dt$$

(A.7)

$$= \frac{1}{\pi |u|^{k+1}} \int_0^{\tau_n+1} |\hat{w}_2^{(k+1)}(t)| \, dt.$$

Compare (A.7) with the lemma; it is sufficient to show that for sufficiently large $n$,

(A.8) $$\int_0^{\tau_n+1} |\hat{w}_2^{(k+1)}(t)| \, dt \leq C,$$

which is equivalent to

(A.9) $$\int_2^{\tau_n+1} |\hat{w}_2^{(k+1)}(t)| \, dt \leq C.$$

We now show (A.9). To do so, we limit our attention to $2 \leq |t| \leq \tau_n + 1$. Recall that $\hat{w}_2(t) = s_2(t) \tilde{w}(t)$, where

$$\tilde{w}(t) = e^{\delta_n t^2/2} \hat{w}_1(t) + \frac{1}{\vartheta_0} \frac{1 - \eta_n}{\eta_n} (e^{\delta_n t^2/2} - 1).$$

First, by the way $\delta_n$ is defined,

(A.10) $$|\tilde{w}(t)| \leq C[|t|^{-\alpha} + t^2 \tau_n^{-(\alpha+2)}] \leq C|t|^{-\alpha}.$$

Second, fix $m = 1, 2, \ldots, k+1$, write

(A.11) $$\tilde{w}^{(m)}(t) = \sum_{j=0}^m (e^{\delta_n t^2/2})^{(m-j)} \hat{w}_1^{(j)}(t) + \frac{1}{\vartheta_0} \frac{1 - \eta_n}{\eta_n} (e^{\delta_n t^2/2})^{(m)}.$$

Recall that $\hat{w}_1(t) = |t|^{-\alpha}$. By elementary calculus, there is a constant $C = C(k) > 0$ such that

(A.12) $$|(e^{\delta_n t^2/2})^{(m)}| \leq C \delta_n t, \qquad |\hat{w}_1^{(m)}(t)| \leq C|t|^{-\alpha}.$$

Combining (A.11) and (A.12) gives

(A.13) $$|\tilde{w}^{(m)}(t)| \leq C \delta_n |t|^{1-\alpha} + C \frac{\delta_n t}{\vartheta_0 \eta_n}$$

$$\leq C \delta_n |t|^{1-\alpha} + C \tau_n^{-(\alpha+1)}, \qquad m = 1, 2, \ldots, k+1.$$

Last, direct calculations show that

(A.14) $$|s_2^{(m)}(t)| \leq C, \qquad m = 0, 1, \ldots, k.$$

Combining (A.10), (A.13) and (A.14) gives

(A.15) $$|\hat{w}_2^{(k+1)}(t)| \leq C[\delta_n |t|^{1-\alpha} + \tau_n^{-(\alpha+1)} + |t|^{-\alpha}]$$

and (A.9) follows directly.

**A.3. Proof of Lemma 6.3.** The first claim follows by the way that $\hat{f}_2$ is constructed. Consider the second claim. Recall that

$$\hat{f}_1(t) = \eta_n e^{-(a^2+1)t^2/2} + e^{-a^2 t^2/2}[(1 - \eta_n) + \vartheta_0 \eta_n \hat{w}_1(t)],$$

and that

$$\hat{f}_2(t) = \eta_n e^{-(a^2+1)t^2/2} + e^{-a_n^2 t^2/2}[(1 - \eta_n) + \vartheta_0 \eta_n \hat{w}_2(t)].$$

Fix $0 \leq m \leq k$. On one hand,

$$|(e^{-a^2 t^2/2})^{(m)}(t)| \leq C|t|^m e^{-a^2 t^2/2}.$$

On the other hand, by the proof of Lemma 6.2,

$$|\hat{w}_1^{(m)}(t)| \leq C, \qquad |\hat{w}_2^{(m)(t)}| \leq C.$$

Combining these gives the claim.

**A.4. Proof of Lemma 6.4.** Write for short $\hat{t} = \hat{t}_n(\gamma)$. By elementary calculus, for any $t > 0$,

(A.16) $$|\varphi_n(t) - 1| \leq \frac{1}{n} \sum_{j=1}^n |e^{itX_j} - 1| \leq \frac{t}{n} \sum_{j=1}^n |X_j|.$$

Note that for sufficiently large $n$, $|\varphi_n(\hat{t}_n)| = n^{-\gamma} \leq 1/2$. Applying (A.16) with $t = \hat{t}_n$ gives

(A.17) $$\hat{t}_n \geq \frac{n}{\sum_{j=1}^n |X_j|} |1 - \varphi_n(\hat{t}_n)| \geq \frac{n/2}{\sum_{j=1}^n |X_j|}.$$

Now, first, by direct calculations and the Hölder inequality,

(A.18) $$|\sigma_0^2(\varphi_n, \hat{t})| \leq \frac{|\operatorname{Re}(\varphi_n(\hat{t})) \operatorname{Re}(\varphi_n'(\hat{t})) + \operatorname{Im}(\varphi_n(\hat{t})) \operatorname{Im}(\varphi_n'(\hat{t}))|}{\hat{t} |\varphi_n(\hat{t})|}$$
$$\leq n^\gamma |\varphi_n'(\hat{t})|/\hat{t},$$

where in the last step we have used $|\varphi_n(\hat{t})| = n^{-\gamma}$. Second, note that for any $t$,

$$(A.19) \qquad |\varphi'(t)| \leq \left|\frac{i}{n}\sum_{j=1}^{n} X_j e^{itX_j}\right| \leq \frac{1}{n}\sum_{j=1}^{n} |X_j|.$$

Combine (A.17), (A.18) and (A.19) and use the Cauchy–Schwarz inequality,

$$(A.20) \qquad |\sigma_0^2(\varphi_n, \hat{t})| \leq 2n^{\gamma}\left(\frac{1}{n}\sum_{j=1}^{n} |X_j|\right)^2 \leq 2n^{\gamma}\left(\frac{1}{n}\sum_{j=1}^{n} X_j^2\right).$$

Hence, to show the claim, it is sufficient to show

$$(A.21) \qquad E\left[\left(\frac{1}{n}\sum_{j=1}^{n} X_j^2\right) \cdot 1_{\{B_n^c\}}\right] \leq C/n.$$

We now show (A.21). Recall that $m_2$ denotes the second moment of $X_1$, we write

$$(A.22) \qquad \frac{1}{n}\sum_{j=1}^{n} X_j^2 = m_2 + \frac{z}{\sqrt{n}},$$

where $z = \sqrt{n}[\frac{1}{n}\sum_{j=1}^{n} X_j^2 - m_2]$. It is seen that $Ez^2 \leq C$, so by the Hölder inequality,

$$(A.23) \qquad \left|E\left[\frac{1}{\sqrt{n}}z \cdot 1_{\{B_n^c\}}\right]\right| \leq \left(\frac{1}{n}Ez^2 \cdot P\{B_n^c\}\right)^{1/2} \leq C/n.$$

Inserting (A.23) into (A.22) gives (A.21). This concludes the proof.

**A.5. Proof of Lemma 6.5.** Before we show the Lemma 6.5, we need some notation and lemmas. Introduce the event

$$(A.24) \qquad D_1 = \left\{W_1(\varphi_n; n) \leq \frac{m_2(\sqrt{(q-2)\log n} + 2m_2)}{\sqrt{n}}\right\},$$

where

$$W_1(\varphi_n; n) = W_1(\varphi_n; n, X_1, X_2, \ldots, X_n) = \sup_{\{|t-t_n| \leq c_0 n^{\gamma-1/2}\}} |\varphi_n'(t) - \varphi'(t)|,$$

$m_2$ is the second moment of $X_1$, and $c_0$ is a constant defined in (6.19). We have the following lemmas.

LEMMA A.1. *Fix $q \geq 4$ and $\gamma \in (0, 1/2)$. For sufficiently large $n$,*

$$P\{D_1^c\} \leq \bar{o}(n^{\gamma-1}).$$

LEMMA A.2. *Fix $q \geq 4$ and $\gamma \in (0, 1/2)$. For sufficiently large $n$,*

(A.25) $$E[|\varphi'_n(\hat{t}_n) - \varphi'(\hat{t}_n)|^2 \cdot 1_{\{D_0 \setminus D_1\}}] \leq C/n.$$

Here, $\bar{o}(n^a)$ denotes a term which equals $o(n^{a+\delta})$ for any $\delta > 0$. The proof of Lemma A.1 is similar to that of Lemma 6.4 of Jin and Cai (2006) so we skip it. Lemma A.2 is the tricky part of the proof of Lemma 6.5 and is proved in Section A.5.1.

We now proceed to prove Lemma 6.5. Write for short $\hat{t}_n = \hat{t}_n(\gamma)$ and $t_n = t_n(\gamma)$. By triangle inequality,

(A.26)
$$\begin{aligned}
E[|\varphi'_n(\hat{t}_n) &- \varphi'(\hat{t}_n)|^2] \\
&\leq E[|\varphi'_n(\hat{t}_n) - \varphi'(\hat{t}_n)|^2 \cdot 1_{\{D_0 \cap D_1\}}] \\
&\quad + E[|\varphi'_n(\hat{t}_n) - \varphi'(\hat{t}_n)|^2 \cdot 1_{\{D_0^c\}}] \\
&\quad + E[|\varphi'_n(\hat{t}_n) - \varphi'(\hat{t}_n)|^2 \cdot 1_{\{D_0 \setminus D_1\}}].
\end{aligned}$$

First, recall that over the event $D_0$ [i.e., (6.19)],
$$|\hat{t}_n - t_n| \leq c_0 n^{\gamma - 1/2},$$
so by the definition of the event $D_1$,
$$|\varphi'_n(\hat{t}_n) - \varphi'(\hat{t}_n)| \leq C\sqrt{\log(n)}/\sqrt{n} \qquad \text{over } D_0 \cap D_1$$
and

(A.27) $$E[|\varphi'_n(\hat{t}_n) - \varphi'(\hat{t}_n)|^2 \cdot 1_{\{D_0 \cap D_1\}}] \leq C \log(n)/n.$$

Second, note that for all $t$,
$$|\varphi'_n(t) - \varphi'(t)| \leq \frac{1}{n} \sum_{j=1}^{n} [|X_j| + m_1] \leq 2m_1 + \frac{1}{n} \sum_{j=1}^{n} (|X_j| - m_1),$$
where $m_1$ is the first moment of $X_1$. It follows that

(A.28)
$$\begin{aligned}
E[|\varphi'_n(\hat{t}_n) &- \varphi'(\hat{t}_n)|^2 \cdot 1_{\{D_0^c\}}] \\
&\leq C\left(E\left[\left(\frac{1}{n}\sum_{j=1}^{n}(|X_j| - m_1)\right)^2 \cdot 1_{\{D_0^c\}}\right] + 2m_1 P\{D_n^c\}\right).
\end{aligned}$$

Moreover, note that $E[\frac{1}{n}\sum_{j=1}^{n}(|X_j| - m_1)]^4 \leq C/n^2$, by the Hölder inequality,
$$E\left[\left(\frac{1}{n}\sum_{j=1}^{n}(|X_j| - m_1)\right)^2 \cdot 1_{\{D_0^c\}}\right]$$

$$\text{(A.29)} \quad \leq \left( E\left[\frac{1}{n}\sum_{j=1}^{n}(|X_j| - m_1)\right]^4 \cdot P\{D_0^c\} \right)^{1/2}$$

$$\leq o(1/n).$$

Combining (A.28) and (A.29) gives

$$\text{(A.30)} \quad E[|\varphi'_n(\hat{t}_n) - \varphi'(\hat{t}_n)|^2 \cdot 1_{\{D_0^c\}}] \leq C/n$$

and the claim follows by inserting (A.25), (A.27) and (A.30) into (A.26).

A.5.1. *Proof of Lemma A.2.* We prove it for the case $\gamma \leq 1/3$ and the case $\gamma > 1/3$ separately.

Consider the case $\gamma \leq 1/3$ first. By the Taylor expansion, for some $\xi$ that falls between $t_n$ and $\hat{t}_n$,

$$\text{(A.31)} \quad \varphi'_n(\hat{t}_n) - \varphi'(\hat{t}_n) = \varphi'_n(t_n) - \varphi'(t_n) + (\varphi''_n(\xi) - \varphi''(\xi)) \cdot (\hat{t}_n - t_n).$$

By direct calculations and the definition of $D_0$,

$$\text{(A.32)} \quad |\varphi''_n(\xi) - \varphi''(\xi)| \leq \frac{1}{n}\sum_{j=1}^{n}(X_j^2 + E[X_j^2]) \leq C \qquad \text{over } D_0.$$

Also, recall that

$$\text{(A.33)} \quad |\hat{t}_n - t_n| \leq c_0 n^{\gamma - 1/2}.$$

Inserting (A.33) and (A.32) into (A.31) gives

$$|\varphi'_n(\hat{t}_n) - \varphi'(\hat{t}_n)| \leq |\varphi'_n(t_n) - \varphi'(t_n)| + Cn^{\gamma - 1/2},$$

which implies

$$|\varphi'_n(\hat{t}_n) - \varphi'(\hat{t}_n)|^2 \leq C(|\varphi'_n(t_n) - \varphi'(t_n)|^2 + n^{2\gamma - 1}).$$

It follows that

$$\text{(A.34)} \quad \begin{aligned} &E[|\varphi'_n(\hat{t}_n) - \varphi'(\hat{t}_n)|^2 \cdot 1_{\{D_0 \setminus D_1\}}] \\ &\qquad \leq C(E[|\varphi'_n(t_n) - \varphi'(t_n)|^2] + n^{2\gamma - 1} \cdot P\{D_0 \setminus D_1\}). \end{aligned}$$

By Lemma A.1 and elementary statistics,

$$\text{(A.35)} \quad P\{D_0 \setminus D_1\} \leq \bar{o}(n^{\gamma - 1}), \qquad E[|\varphi'_n(t_n) - \varphi'(t_n)|^2] \leq C/n,$$

inserting (A.35) into (A.34) gives

$$E[|\varphi'_n(\hat{t}_n) - \varphi'(\hat{t}_n)|^2 \cdot 1_{\{D_0 \setminus D_1\}}] = C/n + \bar{o}(n^{3\gamma - 1})$$

and the claim follows by $\gamma < 1/3$.

Next, consider the case $\gamma \geq 1/3$. Fix $\delta \in (\gamma, 2 - 3\gamma)$ and let
$$K = K(n, c_0, \gamma, \delta) = c_0 n^{\gamma + \delta/2 - 1/2}.$$
Note here that
$$\gamma + \delta/2 - 1/2 > \frac{3\gamma - 1}{2} \geq 0.$$
Lay out a grid $s_k = t_n + (k - K - 1)n^{-\delta/2}$, $k = 1, 2, \ldots, 2K + 1$. Observe that for any $t \in [s_k, s_{k+1}]$,

(A.36)
$$\begin{aligned}|\varphi_n'(t) - \varphi'(t)| &\leq |\varphi_n'(s_k) - \varphi'(s_k)| \\ &\quad + n^{-\delta/2} \cdot \Big(\sup_{|\xi - t_n| \leq c_0 \cdot n^{\gamma - 1/2}} |\varphi_n''(\xi) - \varphi''(\xi)|\Big).\end{aligned}$$

Combining (A.36) with (A.32) gives
$$|\varphi_n'(t) - \varphi'(t)| \leq |\varphi_n'(s_k) - \varphi'(s_k)| + Cn^{-\delta/2} \quad \text{over } D_0.$$
Now, note that the endpoints of the grid are
$$t_n \pm Kn^{-\delta/2} = t_n \pm c_0 n^{\gamma - 1/2}$$
and that over the event $D_0$,
$$|\hat{t}_n - t_n| \leq c_0 n^{\gamma - 1/2};$$
it follows that
$$|\varphi_n'(\hat{t}_n) - \varphi'(\hat{t}_n)| \leq \max_{\{1 \leq k \leq 2K+1\}} |\varphi_n'(s_k) - \varphi'(s_k)| + Cn^{-\delta/2}.$$
Therefore, by the Cauchy–Schwarz inequality,

(A.37) $|\varphi_n'(\hat{t}_n) - \varphi'(\hat{t}_n)|^2 \leq C\Big(\Big(\max_{\{1 \leq k \leq 2K+1\}} |\varphi_n'(s_k) - \varphi'(s_k)|^2\Big) + n^{-\delta}\Big).$

Recall that

(A.38) $$P\{D_0 \setminus D_1\} \leq \bar{o}(n^{\gamma - 1}).$$

It follows from (A.37) and (A.38) that

(A.39)
$$\begin{aligned}E[|\varphi_n'(\hat{t}_n) &- \varphi'(\hat{t}_n)|^2 \cdot 1_{\{D_0 \setminus D_1\}}] \\ &\leq C\Big(E\Big[\Big(\max_{\{1 \leq k \leq 2K+1\}} |\varphi_n'(s_k) - \varphi'(s_k)|^2\Big) \cdot 1_{\{D_0 \setminus D_1\}}\Big] \\ &\quad + n^{-\delta} P\{D_0 \setminus D_1\}\Big) \\ &= C \sum_{k=1}^{2K+1} E[|\varphi_n'(s_k) - \varphi'(s_k)|^2 \cdot 1_{\{D_0 \setminus D_1\}}] + o(n^{-1}),\end{aligned}$$

where in the last step we have used $\delta > \gamma$.

Now, for any $k = 1, 2, \ldots, 2K+1$, observe that by elementary statistics,
$$E[|\varphi_n'(s_k) - \varphi'(s_k)|^4] \leq C/n^2.$$

By the Hölder inequality and (A.38),
$$\begin{aligned} E[|\varphi_n'(s_k) - \varphi'(s_k)|^2 \cdot 1_{\{D_0 \setminus D_1\}}] \\ \leq (E|\varphi_n'(s_k) - \varphi'(s_k)|^4 \cdot P\{D_0 \setminus D_1\})^{1/2} \\ \leq \bar{o}(n^{(\gamma-3)/2}), \end{aligned}$$

so by $K \leq C n^{\gamma + \delta/2 - 1/2}$

$$\sum_{k=1}^{K} E[|\varphi_n'(s_k) - \varphi'(s_k)|^2 \cdot 1_{\{D_0 \setminus D_1\}}] \leq \bar{o}(Kn^{(\gamma-3)/2})$$
(A.40)
$$= \bar{o}(n^{3\gamma/2 + \delta/2 - 2}).$$

Recall $\delta < 2 - 3\gamma$, it follows from (A.40) that

(A.41) $$\sum_{k=1}^{K} E[|\varphi_n'(s_k) - \varphi'(s_k)|^2 \cdot 1_{\{D_0 \setminus D_1\}}] = o(1/n)$$

and the claim follows by plugging (A.41) into (A.39).

**A.6. Proof of Lemma 6.9.** Consider the first claim. Write for short $\bar{t}_n = \sqrt{2\gamma \log n}$. By the definition and elementary Fourier analysis,

(A.42)
$$\begin{aligned} \tilde{\varphi}^R(t) = (1 - \varepsilon_n) e^{-1/2(\sigma_0/\bar{\sigma}_0)^2 t^2} \cos\left(t \frac{u_0 - \bar{u}_0}{\bar{\sigma}_0}\right) \\ + \varepsilon_n \int e^{-1/2(\sigma/\bar{\sigma}_0)^2 t^2} \cos\left(t \frac{u - \bar{u}_0}{\bar{\sigma}_0}\right) h(u|\sigma) \, dH(\sigma). \end{aligned}$$

By Lemma 6.6, we have that over the event $B_n$,

(A.43) $\quad |\hat{\sigma}_0 - \bar{\sigma}_0| \leq C \log^{1/2}(n) n^{\gamma - 1/2}, \qquad |\hat{u}_0 - \bar{u}_0| \leq C \log(n) n^{\gamma - 1/2}.$

As a result, by the Taylor expansion and that $\tilde{t}_n = \frac{\bar{\sigma}_0}{\hat{\sigma}_0} \bar{t}_n$,

(A.44)
$$\begin{aligned} |\tilde{\varphi}^R(\tilde{t}_n) - \tilde{\varphi}^R(\bar{t}_n)| \lesssim |(\tilde{\varphi}^R)'(\bar{t}_n)| \cdot |\tilde{t}_n - \bar{t}_n| \\ \leq C \log(n) n^{\gamma - 1/2} |(\tilde{\varphi}^R)'(\bar{t}_n)|, \end{aligned}$$

where

(A.45) $\quad |(\tilde{\varphi}^R)'(\bar{t}_n)| \leq C \bar{t}_n e^{-1/2(\sigma_0/\bar{\sigma}_0)^2 \cdot \bar{t}_n^2} \leq C \log^{1/2}(n) n^{-\gamma}.$

Combining (A.44) and (A.45) gives the first claim.

Consider the second claim. Introduce a bridging quantity

$$E\left[\cos\left(\bar{t}_n \frac{X_1 - u_0}{\sigma_0}\right)\right]. \tag{A.46}$$

By the triangle inequality,

$$|(1 - \varepsilon_n) - n^\gamma \tilde{\varphi}^R(\bar{t}_n)| \leq \mathrm{I} + \mathrm{II}, \tag{A.47}$$

where $\mathrm{I} = |(1 - \varepsilon_n) - n^\gamma E[\cos(\bar{t}_n \frac{X_1 - u_0}{\sigma_0})]|$ and $\mathrm{II} = n^\gamma |E[\cos(\bar{t}_n \frac{X_1 - u_0}{\sigma_0})] - \tilde{\varphi}^R(\bar{t}_n)|$. Consider I first. By direct calculations and $e^{\bar{t}_n^2/2} = n^\gamma$,

$$\begin{aligned}(1 - \varepsilon_n) - n^\gamma E\left[\cos\left(\bar{t}_n \frac{X_1 - u_0}{\sigma_0}\right)\right] \\ = -\varepsilon_n \int e^{-1/2[(\sigma/\sigma_0)^2 - 1]\bar{t}_n^2}\left[\int \cos\left(\bar{t}_n \frac{u - u_0}{\sigma_0}\right) h(u|\sigma)\right] dH(\sigma).\end{aligned} \tag{A.48}$$

Note that by elementary Fourier analysis,

$$\int \cos\left(t\frac{u - u_0}{\sigma_0}\right) h(u|\sigma) = \mathrm{Re}\left(\tilde{h}\left(\frac{t}{\sigma_0}\bigg|\sigma\right)\right).$$

Since $H$ is eligible and obeys the constraint (2.3), we have

$$\begin{aligned}\left|(1 - \varepsilon_n) - n^\gamma E\left[\cos\left(\bar{t}_n \frac{X_1 - u_0}{\sigma_0}\right)\right]\right| \\ \leq \varepsilon_n \int e^{-1/2[(\sigma/\sigma_0)^2 - 1]\bar{t}_n^2}\left|\tilde{h}\left(\frac{\bar{t}_n}{\sigma_0}\bigg|\sigma\right)\right| dH(\sigma) \\ \leq C\varepsilon_n \bar{t}_n^{-\alpha}.\end{aligned} \tag{A.49}$$

Consider II next. It follows from the proof of Theorem 2.2 [i.e., (6.18)] that

$$\begin{aligned}|\bar{\sigma}_0 - \sigma_0| &\leq C\varepsilon_n \log^{-(\alpha+2)/2}(n), \\ |\bar{u}_0 - u_0| &\leq C\varepsilon_n \log^{-(\alpha+1)/2}(n).\end{aligned} \tag{A.50}$$

Compare (A.48) with (A.42),

$$\begin{aligned}\left|\tilde{\varphi}^R(\bar{t}_n) - E\left[\cos\left(\bar{t}_n \frac{X_1 - u_0}{\sigma_0}\right)\right]\right| \\ \leq Cn^{-\gamma}\bigg[(1 - \varepsilon_n)(|\bar{\sigma}_0^2 - \sigma_0^2|\bar{t}_n^2 + |\bar{u}_0 - u_0|\bar{t}_n) \\ + \varepsilon_n \int (\sigma^2 \bar{t}_n^2|\bar{\sigma}_0^2 - \sigma_0^2| + |u\bar{t}_n||\bar{u}_0 - u_0|) \, dH(u, \sigma)\bigg].\end{aligned}$$

Note that $E|u| \leq C$ and $E|\sigma^2 - \sigma_0^2| \leq C$, it follows from (A.50) that

$$(A.51) \qquad \left| \tilde{\varphi}^R(\bar{t}_n) - E\left[\cos\left(\bar{t}_n \frac{X_1 - u_0}{\sigma_0}\right)\right] \right| \leq C\varepsilon_n n^{-\gamma} \log^{-\alpha/2}(n).$$

Inserting (A.49) and (A.51) to (A.47) gives

$$|(1-\varepsilon_n) - n^\gamma \tilde{\varphi}^R(\bar{t}_n)| \leq C\varepsilon_n \log^{-\alpha/2}(n).$$

This concludes the proof of the second claim of the lemma.

**Acknowledgments.** We would like to thank Wenguang Sun for his help on programming and numerical simulations. We also thank the Associate Editor and two referees for thorough and useful comments which have helped to improve the presentation of the paper.

## REFERENCES

ABRAMOVICH, F., BENJAMINI, Y., DONOHO, D. and JOHNSTONE, I. (2006). Adapting to unknown sparsity by controlling the false discovery rate. *Ann. Statist.* **34** 584–653. MR2281879

BENJAMINI, Y. and HOCHBERG, Y. (1995). Controlling the false discovery rate: A practical and powerful approach to multiple testing. *J. Roy. Statist. Soc. Ser. B* **57** 289–300. MR1325392

BENJAMINI, Y. and HOCHBERG, Y. (2000). On the adaptive control of the false discovery rate in multiple testing with independent statistics. *J. Educ. Behav. Statist.* **25** 60–83.

BENJAMINI, Y., KRIEGER, A. M. and YEKUTIELI, D. (2006). Adaptive linear step-up procedures that controls the false discovery rate. *Biometrika* **93** 491–507. MR2261438

BLANCHARD, G. and ROQUAIN, É. (2007). Adaptive FDR control under independence and dependence. Available at arxiv:0707.0536v2.

CAI, T., JIN, J. and LOW, M. (2007). Estimation and confidence sets for sparse normal mixtures. *Ann. Statist.* **35** 2421–2449. MR2382653

CELISSE, A. and ROBIN, S. (2008). A leave-p-out based estimation of the proportion of null hypotheses. Available at arxiv:0804.1189. MR2411944

DONOHO, D. and JIN, J. (2004). Higher criticism for detecting sparse heterogeneous mixtures. *Ann. Statist.* **32** 962–994. MR2065195

DONOHO, D. and JIN, J. (2006). Asymptotic minimaxity of false discovery rate thresholding for sparse exponential data. *Ann. Statist.* **34** 2980–3018. MR2329475

DONOHO, D. and LIU, R. C. (1991). Geometrizing rates of convergence, II. *Ann. Statist.* **19** 633–667. MR1105839

EFRON, B. (2004). Large-scale simultaneous hypothesis testing: The choice of a null hypothesis. *J. Amer. Statist. Assoc.* **99** 96–104. MR2054289

EFRON, B. (2008). Microarrays, empirical Bayes and the two-groups model. *Statist. Sci.* **23** 1–22. MR2431866

EFRON, B., TIBSHIRANI, R., STOREY, J. and TUSHER, V. (2001). Empirical Bayes analysis of a microarray experiment. *J. Amer. Statist. Assoc.* **96** 1151–1160. MR1946571

ERDELYI, A. (1956). *Asymptotic Expansions*. Dover, New York. MR0078494

FAN, J. (1991). On the optimal rates of convergence for nonparametric deconvolution problems. *Ann. Statist.* **19** 1257–1272. MR1126324

Finner, H., Dickhaus, T. and Roters, M. (2009). On the false discovery rate and on asymptotically optimal rejection curve. *Ann. Statist.* **37** 596–618. MR2502644

Genovese, C. and Wasserman, L. (2004). A stochastic process approach to false discovery control. *Ann. Statist.* **32** 1035–1061. MR2065197

Ibragimov, I. A., Nemirovskii, A. S. and Khas'minskii, R. Z. (1986). Some problems on nonparametric estimation in Gaussian white noise. *Theory Probab. Appl.* **31** 391–406. MR0866866

Jin, J. and Cai, T. (2006). Estimating the null and the proportion of nonnull effects in large-scale multiple comparisons. Available at arxiv.math/0611108v1. MR2325113

Jin, J. and Cai, T. (2007). Estimating the null and the proportion of nonnull effects in large-scale multiple comparisons. *J. Amer. Statist. Assoc.* **102** 495–506. MR2325113

Jin, J. (2008). Proportion of nonzero normal means: Universal oracle equivalences and uniformly consistent estimations. *J. R. Stat. Soc. Ser. B Stat. Methodol.* **70**(3) 461–493.

Mallat, S. (1998). *A Wavelet Tour of Signal Processing*, 2nd ed. Academic Press, New York. MR1614527

Meinshausen, M. and Rice, J. (2006). Estimating the proportion of false null hypothesis among a large number of independent tested hypotheses. *Ann. Statist.* **34** 373–393. MR2275246

Müller, P., Parmigiani, G., Robert, C. and Rousseau, J. (2004). Optimal sample size for multiple testing: The case of gene expression microarrays. *J. Amer. Statist. Assoc.* **99** 990–1001. MR2109489

Newton, M., Kendziorski, C., Richmond, C., Blattner, F. and Tsui, K. (2001). On differential variability of expression ratios: Improving statistical inference about gene expression changes from microarray data. *J. Comput. Biol.* **8** 37–52.

Neuvial, P. (2008). Asymptotic properties of false discovery rate controlling procedures under independence. *Electron. J. Stat.* **2** 1065–1110. MR2460858

Silverman, B. W. (1986). *Density Estimation for Statistics and Data Analysis*. Chapman and Hall, London. MR0848134

Storey, J. D. (2002). A direct approach to false discovery rate. *J. R. Stat. Soc. Ser. B Stat. Methodol.* **64** 479–498. MR1924302

Sun, W. and Cai, T. (2007). The oracle and compound decision rules for false discovery rate control. *J. Amer. Statist. Assoc.* **102** 901–912. MR2411657

Swanepoel, J. W. H. (1999). The limiting behavior of a modified maximal symmetric 2s-spacing with applications. *Ann. Statist.* **27** 24–35. MR1701099

van der Laan, M., Dudoit, S. and Pollard, K. (2004). Multiple testing (III): Augmentation procedures for control of the generalized family-wise error rate and tail probabilities for the proportion of false positives. Technical report, Dept. Biostatistics, Univ. California, Berkeley.

West, M. (1987). On scale mixtures of normal distributions. *Biometrika* **3** 646–648. MR0909372

Zhang, C.-H. (1990). Fourier methods for estimating mixing densities and distributions. *Ann. Statist.* **18** 806–831. MR1056338

Department of Statistics
The Wharton School
University of Pennsylvania
Philadelphia, Pennsylvania 19104
USA
E-mail: tcai@wharton.upenn.edu

Department of Statistics
Carnegie Mellon University
Pittsburgh, Pennsylvania 15213
USA
E-mail: jiashun@stat.cmu.edu